\documentclass[12pt]{amsart}
\usepackage{amsmath,amssymb,epsfig,amssymb,amsthm,epic,eepic,color,amsfonts}
\usepackage{hyperref}
\usepackage[all]{xy}

\usepackage{fontenc}  
\usepackage[utf8]{inputenc}
\usepackage{enumerate}

\newtheorem{thm}{Theorem}[section]

\newtheorem{prop}[thm]{Proposition}
\newtheorem{defn}[thm]{Definition}
\newtheorem{lemme}[thm]{Lemma}

\newtheorem{important-rem}[thm]{Important remark}

\newtheorem{rien}[thm]{}
\newtheorem*{defnt}{Definition}

\numberwithin{equation}{section}

\newcommand{\be}{\begin{enumerate}}
\newcommand{\ee}{\end{enumerate}}
\newcommand{\bi}{\begin{itemize}}
\newcommand{\ei}{\end{itemize}}

\def\R{\mathbb{R}}

\def\Z{\mathbb{Z}}

\def\F{\mathcal{F}}
\def\D{\mathcal{D}}
\def\M{\mathcal{M}}

\def\L{\mathcal{L}}
\def\A{\mathbb A}

\def\om{\omega}

\def\ga{\gamma}    
\def\Ga{\Gamma}

\def\al{\alpha}
\def\be{\beta}

\def\De{\Delta}
\def\vp{\varphi}

\def\la{\lambda}

\def\si{\sigma}
\def\Si{\Sigma}

\def\ep{\varepsilon}

\def\ds{\displaystyle}
\def\p{\partial}

\def\nd{\noindent}
\def\bull{\hfill$\Box$\\}
\def\proof{\nd {\bf Proof.\ }}

\begin{document}
\vskip 1cm
\begin{center}
 About the diffeomorphisms of the 3-sphere and a famous theorem of Cerf ($\Ga_4=0$) 
 \vskip 1cm
 Fran\c cois Laudenbach
\end{center}

\title{}
\author{}
\address{Laboratoire de
math\'ematiques Jean Leray,  UMR 6629 du
CNRS, Facult\'e des Sciences et Techniques,
Universit\'e de Nantes, 2, rue de la
Houssini\`ere, F-44322 Nantes cedex 3,
France.}
\email{ francois.laudenbach@univ-nantes.fr}



\begin{abstract}  Using a rigidity property of the foliations of $S^2\times [0,1]$
 that are defined by a non-vanishing closed 
one-form, we give a rather simple proof of a theorem due J. Cerf, going back to 1968, that the group of direct diffeomorphisms of  $S^3$ is connected.  
\end{abstract}
\maketitle
\thispagestyle{empty}
${}$ \hfill


The famous theorem  of Jean Cerf in question in the title states the following.\\

\emph{Every diffeomorphism of the 3-sphere preserving the orientation is isotopic to the identity}.\\

It follows that every diffeomorphism of $S^3$ extends to the 4-ball. Then, $\Ga_4:={\rm Diff}S^3/\rho({\rm Diff} D^4)$
is a trivial group; here, $\rho$ stands for the  restriction morphism to the boundary of $D^4$. So, $\Ga_4=0$ is a short name 
for the real theorem of Cerf.
Actually in 1992, using recent results at the time in 3-dimensional contact 
geometry and the theory of holomorphic discs,  Y. Eliashberg gave a direct proof of $\Ga_4= 0$  
\cite[section 6]{eliash},
avoiding Cerf's theorem.  

In 1979, I solved in \cite{lauden}, with the help of Samuel Blank, a problem raised by J. Moser \cite{moser}: 
\begin{equation*} 
(\rm A)\ \left\{
\begin{array}{l}
\text{\it On a compact 3-dimensional manifold, two non-vanishing closed one-forms}\\
 \text{\it that are tangent to the boundary and cohomologous are isotopic.}
 \end{array}\right.
\end{equation*}

Our proof did not use Cerf's theorem.
Applying (A) to $S^2\times [0,1]$ 
 immediately implies the original theorem of Cerf, not just $\Ga_4=0$.\footnote{ For completeness, one should apply
some classical fibration theorems also due to Cerf \cite[Appendice]{cerf}.} Nevertheless, I should confess that this 
isotopy theorem of 1-forms is somehow technical. A simpler proof of (A) is given by N.V. Qu\^e and R. Roussarie \cite{roussarie},
but depending on Cerf's theorem. 

We name Theorem (A') the particular case of Theorem (A) where the ambient manifold
is  $S^2\times [0,1]$. The present paper aims to give a less technical proof of Theorem (A') than \cite{lauden}, 
and hence of Cerf's Theorem.\\

\nd {\bf Acknowledgement.} I am grateful to Robert Roussarie for his reading of a preliminary version of this text,
 his advice about the saddles of type $X$ and other comments. 
 I thank the project COSY ANR-21-CE40-0002 for its support.

\section{Preliminaries and plan of proof}

In the rest of the paper,  the ambient manifold will be $S^2\times [0,1]$ and 
$z$ will denote the variable in its second factor.

\begin{rien}{\rm 
  Every  non-vanishing closed one-form $\om$ on $S^2\times [0,1]$ that is tangent to the boundary reads 
   $\vp^*dz$ for some diffeomorphism $\vp$ that is equal to identity on $S^2\times \{0\}$. 
 This 
 follows from the fact that 
 every orbit of a vector field  $X$ that everywhere fulfills $\om(X)>0$
 goes from $S^2\times\{0\}$ to $S^2\times\{1\}$.
Moreover, if $\vp^*dz= dz$ then $\vp$ is isotopic to the identity
among the diffeomorphisms of the same type.\footnote{ These facts hold true on $S^n\times[0,1]$ 
for every dimension. }  
}
\end{rien}

\begin{rien}{\sc  Moser's method. }\label{moser-trick}
{\rm Given two non-vanishing closed 1-forms $\om_0$ and $\om_1$
on  $S^2\times [0,1]$, tangent to the boundary, a point $a\in S^2\times [0,1]$ 
is said to be a \emph{negative contact}  (resp. a \emph{positive contact}) 
if there exists $\mu>0$ (resp. $\mu<0$) such that $\om_0(a)\pm\mu\,\om_1(a)=0$. In other words,
the kernels of $\om_0(a)$ and $\om_1(a)$ coincide but their co-orientations are opposite (resp. equal).
The locus of negative (resp. positive) contacts will be denoted by $C_-$ (resp. $C_+$.)}
\end{rien}

\begin{lemme}\label{moser-trick2}
{If the two non-vanishing closed one-forms are cohomologous and their mutual negative contact 
locus $C_-$ is empty, 
then the two forms are isotopic.}
\end{lemme}

 Indeed, the barycentric combination
$\om_t:= (1-t)\om_0+t\om_1$ provides
a path $(\om_t)_{t\in [0,1]}$ of  closed 1-forms, from $\om_0$ to $\om_1$, that are cohomologous, tangent to the 
boundary and nowhere singular. In that case,
Moser's {\it isotopy theorem} applies, as in the case of volume forms \cite{moser}.
Namely,  there is an autonomous flow $(\psi_t)_{t\in[0,1]}$ 
that conjugates $\om_t$ to $\om_0$ for every $t\in[0,1]$.\footnote{ In the case of non-vanishing closed one-forms, 
a very simple proof is given in \cite[Appendice I]{lauden}.} 

Therefore, our method to prove the isotopy Theorem (A') on $S^2\times [0,1]$
will be to cancel the locus $C_-$ of negative contacts  between $\om_0:=dz$ 
and $\om_1:=\vp^*dz$.\\

\begin{rien} {\sc Reduction to a more special case.}
{\rm Let $f$ be the function which is equal to the projection onto the factor $[0,1]$;
we set $\ell:= \vp^*f= f\circ \vp$. 
From now on, the two forms  $\om_0$ and $\om_1$ will  be respectively denoted by $\om_\F$ and $\om_\L$; 
and the foliations they define will be respectively named $\F$ and $\L$. The leaves of $\om_\F$ 
 will be named the \emph{level sets} (understood
   of $f$)  while the level sets of
$\ell$ will continue to be named the {\it leaves} (understood 
of $\L$.)

Without loss of generality, to prove Theorem (A') we may assume the following in the rest of the paper:
\begin{equation}\label{no-saddle}
\left\{\begin{array}{l}
 \text{The restriction of }f
\text{ to every leaf of }\ell, \text{ close enough to the boundary }S^2\times\{0,1\}, \\
\text{has only two critical points, a maximum and a minimum, and these belong to }C_+.\\
\end{array} \right.
\end{equation}
 }
\end{rien}
\medskip

\begin{rien} {\sc Plan of the proof of Theorem (A').} \label{plan}
{\rm Generically,
the contact points (positive or negative) have a Morse index: 0, 1 or 2. The contact points of index 1 are named {\it saddles.}
They have also a \emph{type}: $Y, \ \la, \text{ or }X$. A saddle  $s$ on a leaf  $L\in\L$ is said to be 
of type $\la$ (resp. $Y$) if the connected component of $L\cap [f(s)-\ep, f(s)+\ep]$  which contains $s$ looks like a \emph{pair-of-pants} (resp. a reversed one.) 
A saddle is said to be of type $X$ if it is 
 the common limit of two sequences of saddles, 
one of  $\la$-saddles
and the other of $Y$-saddles, the pair $(f,\ell)$ being kept fixed (Definition \ref{saddle--text}.) 
This necessarily involves a pair of saddles both of type $X$ in the same connected component of level curves 
(Lemma \ref{X-pair}.)
One will speak of 
$\la$-, $Y$-, and $X$-saddle respectively.

Section \ref{basics} is devoted to basics on generic properties. Section \ref{elementary} reviews 
elementary isotopies whose effect is to simplify
$C_-$. One of them will allow us to kill all the finitely many $X$-saddles. Typically, this phenomenon 
occurs
in the setting {\it saddle-center-saddle} (Definition \ref{scs-conf}.)

The next operation 
will consist of pushing the saddles of type $Y$ to levels higher than
 all saddles of type  $\la$. 
Unfortunately,
this operation is obstrued by connecting orbits $Y\to\la$ of the $\L$-gradient of $f$, a phenomenon
that generically appears finitely many times. Here, the $\L$-gradient of $f$ is  the orthogonal projection of 
$\nabla f=\p_z$  
onto the leaves of $\ell$ with respect to an understood  Riemaniann metric; it will be denoted by $\nabla_\L f$. 

These obstructions are destroyed
thanks to the main operation of the paper, namely the \emph{Dehn modification}. Section \ref{dehn} will be devoted to this
type of 3-dimensional surgery considered in our specific setting. 
Section \ref{main-proof} makes use of all these tools to prove the main result, namely Theorem (A').
}
\end{rien}
\section{Basics on pairs of non-vanishing closed one-forms}\label{basics}

\begin{rien} {\sc Regularity.} \label{regul}{\rm 
We set $C=C_-\cup C_+$. 
With local coordinates $(x,y)$ on $S^2$, the contact locus is defined by the following system: 
\begin{equation}\label{reg}
C =\left\{ \frac{\p\ell}{\p x}(x,y,z)=0, \ \frac{\p\ell}{\p y}(x,y,z)=0 \right\}.
\end{equation} 
The subsets $C_+$ (resp. $C_-$) are defined by adding the inequation $\ds\frac{\p\ell}{\p z}(x,y,z)>0\ (\text{resp.}<0)$.
By the assumption (\ref{no-saddle}), $C_-$ does not approach  the boundary and hence is compact.

By the transversality theorem of  
Thom in jet spaces, some approximation of $\vp$ makes maximal the rank of
the linearized system  associated with system (\ref{reg})
at every point of $S^2\times ]0,1[$. In this case 
 the pair $(\om_\F, \om_\L)$---or $(f,\ell)$---is said to be {\it regular} and $C$ is a smooth curve.
The open subset of points in $C$ where 
\begin{equation}
{\ds\De:= \frac{\p^2\ell}{\p x^2}\frac{\p^2\ell}{\p y^2}-\left( \frac{\p^2\ell}{\p x\p y}\right)^{\! 2}}
\end{equation}
 is non-zero
is the locus where $C$ is transverse to both foliations $\F$ and $\L$; overthere, their contact is quadratic. If 
$a$ is such a point and $L$ is the leaf of $a$, this point has a Morse index in $\{0,1,2\}$ as a critical point 
of the function $f_{L}$, that is the restriction of $f$ to $L$. The critical points of $f_L$ are named \emph{minimum, saddle,
maximum,} depending on their Morse index. 

The remaining points of $C$ are called {\it inflection points}. Generically, the equation $\De(x,y,z)=0$ is regular
and hence, the inflections are isolated in $C$ and do not approach the boundary. So, there are finitely many of them.
 At an inflection,  the tangency of both $\F$ and $\L$ with $C$ is quadratic 
 while the mutual tangency of $\F$ and $\L$ is cubic.  
 
 By the {\it versal unfolding} theory \cite{mather} 
\cite{po}, there are coordinates in a neighborhood of an inflection  where the two functions $f$ and $\ell$ read
\begin{equation}\label{inflection-f}
f(x,y,z)= z \text{ and }\ell(x,y,z)= x^3\pm y^2+\la(z)\,x  +\mu(z).
\end{equation}
}
\end{rien}
\medskip
The sign in the above formula depends on the Morse indices of contact points nearby. The Hessian of $f_L$
at an inflection $I\in L$ has index 0 (resp. 1), and is said to be a \emph{saddle-min} (resp. \emph{saddle-max})
 inflection. If the restriction $f_C$ of $f$ to $C$ is locally minimal (resp. maximal) at $p\in C$, $p$ is said to be
 a \emph{birth} (resp. \emph{cancellation}) inflection.

Since $C_-$ is generically a closed---in general non-connex---curve, 
 every  of its connected component carries at least two inflections, one maximum and one minimum of the restriction
 of $f$ to $C_-$.
 
The Morse index of a quadratic contact is constant on each arc of $C$ ending at an inflection 
or boundary points.
When two arcs in $C$ of quadratic contact have a common inflection 
in their closure, 
their indices differ by 1. So, one of these two arcs
has a constant index 1; 
every contact point of this arc is a saddle,
and the other arc has the index of an $f_L$-extremum.\\
    
 \begin{rien}{\sc Excellence.} {\rm The pair of functions $(f,\ell)$ is said to be \emph{excellent} if the following 
 conditions are fufilled.
 \begin{enumerate}
 \item [(E1)]  Let $a$ and $b$ be two distinct quadratic contact points at the same  level and in the same leaf;
 let $(a(t), b(t))$ is a pair of local parametrizations of the respective contact arcs such that 
 $a(0)=a, b(0)=b \text{ and } \ell(a(t))=\ell(b(t))$ for every $t$ close to $0$.
 Then we have $\ds\frac {d }{dt}[f(a(t))-f(b(t))]\neq 0 $ at $t=0$. 
 \smallskip
 \item[(E2)] There are neither three contact points in the same leaf and in the same level set
 nor two pairs of contact points
 in the same leaf at two distinct levels.
 \item[(E3)] Let $\si$ be an inflection point. Then the leaf of $\si$ (resp.  its level set) contains 
  no contact point at the same level (resp. in the same leaf) as $\si$.
 \end{enumerate}
  }
 \end{rien}
 
 \begin{lemme} In the space of smooth functions $\ell: S^2\times[0,1]\to [0,1] $
 such that the pair $(f,\ell)$ fulfils the assumption {\rm (\ref{no-saddle})}, the subspace 
 such that the pair $(f,\ell)$ is regular and excellent is an open dense set. 
 \end{lemme} 
 
 \proof Such a statement is well known to the topologists when one speaks of a generic path of functions,
 the variable $z$ being the time. But every path of functions is not the primitive of a  non-vanishing one-form $\om_\L$
 tangent to the boundary. 
 The main difference comes from the points where 
 $\ds\frac {\p \ell}{\p z}=0$. 
 In our setting, $\ds\frac {\p \ell}{\p z}$ is not vanishing at every point where the two other partial 
 derivatives vanish. That means that we work in an open set of the general space of smooth real functions on $S^2$. So,
 Thom's transversality theorems
 in bi-one-jet\footnote{ A bi-one-jet is just a pair of two one-jets.}  spaces of real functions apply. \bull
 
 \section{Elementary isotopies and application.}\label{elementary}
 
 Here is the list of the isotopies of $\ell$ that we consider as \emph{elementary}: 
  cancellation---or creation---of a \emph{simple loop} in the contact locus $C$;
 bypass of the \emph{cusp} singularity $x^4$;
 exchange (that is another way of bypassing a cusp.) 
 To this list, we add the \emph{isotopies along a satured set}, though they are not related to singularity theory.
 
 As announced, these techniques will allow us to cancel the saddles of type $X$ (in the sense of subsection \ref{plan}.) Before starting
 with the description of these isotopies, some notation is needed.
 
 \begin{rien} \label{notation-ball}
 {\rm Given a leaf $L\in\L$ and an extremum $m$ of $f_{L}$, 
 the \emph{cone} of $m$, denoted by $C_\L(m)$, is the largest open disc in $L$---if it exists---that contains 
 $m$ and is bounded by a \emph{singular} simple closed curve $\p C_\L(m)$ in a level set of $f$ such that:
 
 (i) $\p C_\L(m)$ contains a saddle (or an inflection) and bounds a disc $\De(m)$ in its level set;
 
 (ii) $\De(m)\cup C_\L(m)$ bounds a (topological) ball $B_\L(m)$ in the ambient manifold.
 } 
 \end{rien}
 
 Note that if $L$ is close to the boundary of $S^2\times [0,1]$, by (\ref{no-saddle}) the cone of the minimum 
 (resp. maximum) is not 
 defined since the leaf contains no saddles.\\
 
 \begin{defn}\label{def-simple} {\sc(Simple loop.)} A simple closed curve 
 $\Ga$ in 
 the contact locus $C$ of the pair $(f,\ell)$ of smooth 
 functions on $S^2\times [0,1]$ is said to be a \emph{simple loop} if the following conditions are fulfilled.
 \begin{enumerate}[{\rm (1)}]
 \item $\Ga$ is made of two open arcs, $\al$ and $\beta$, of quadratic contact points of the same signs
 and two inflection points, $I_0$ and $I_1$, which make the topological closure of each of the arcs $\al$ and $\beta$;
 say $\al$ is of Morse index $1$; the index of $\beta$, $0 \text{ or }2$, depends on the sign in front of $y^2$ in 
 formula {\rm(\ref{inflection-f})}.
 \item Every leaf $L$ that crosses $\al$ at a saddle $s_L$\footnote{ By the Stokes formula this point is unique when it exists.} also crosses $\beta$ at a unique extremum $m_L$ 
 and conversely. Moreover, $s_L$ 
 belongs to $\p C_\L(m_L)$.
 \item There is a unique $\L$-gradient line\footnote{ See the definition in subsection \ref{plan}.} $\ga_L$ joining $s_L$ to $m_L$ on every leaf $L$ meeting $\al\cup\beta$.
 \item The only arc of definite quadratic contacts located in the ball $B_\L(m_L)$ is $\beta\cap B_\L(m_L)$.
 \end{enumerate}
  \end{defn}
  By formula (\ref{inflection-f}), if $s_L$ goes to $I_j$ ($j=0 \text{ or }1$) then $\ga_L$ with its field of unit tangent vectors
  goes to $I_j$ with $\pm \p_x$. Note that, by permuting $f \text{ and }\ell$, there are $\F$-gradient lines connecting pair of points in 
  $\al\times\beta$ lying in the same level set of $f$.

  \begin{prop}\label{cancel-loop} Let $\Ga$ be a \emph{simple loop} in the contact locus of the pair $(f,\ell)$. 
  Then, there exists 
  a smooth closed 3-ball $N$ which is an arbitrary small neighbourhood of the union $\bigcup_{m_L\in \beta} B_\L(m_L)$
  satisfying the following two properties.
  \begin{enumerate}[{\rm (1)}]
  \item The boundary is foliated by smooth closed curves 
 drawn by $\L\cap \p N$, 
  with exactly two singularities 
  at $\max\ell_{\p N}$ and $\min\ell_{\p N}$.\footnote{ These two points are close to $I_0$ and $I_1$
  respectively, depending
  on the distance of $\p N$ to $\cup_{m_L\in \beta} B_\L(m_L)$.}
  
  \item The germ of $\ell$ along $\p N$ extends to $N$ as a  smooth function $\ell'$ without contact point with $f$ in $N$.
  Moreover, $\ell'$ is isotopic to $\ell$ by an isotopy supported in the interior of $N$.
   \end{enumerate}
  \end{prop}
  
   \begin{center}
 \begin{figure}[h]
\includegraphics[scale =.7]{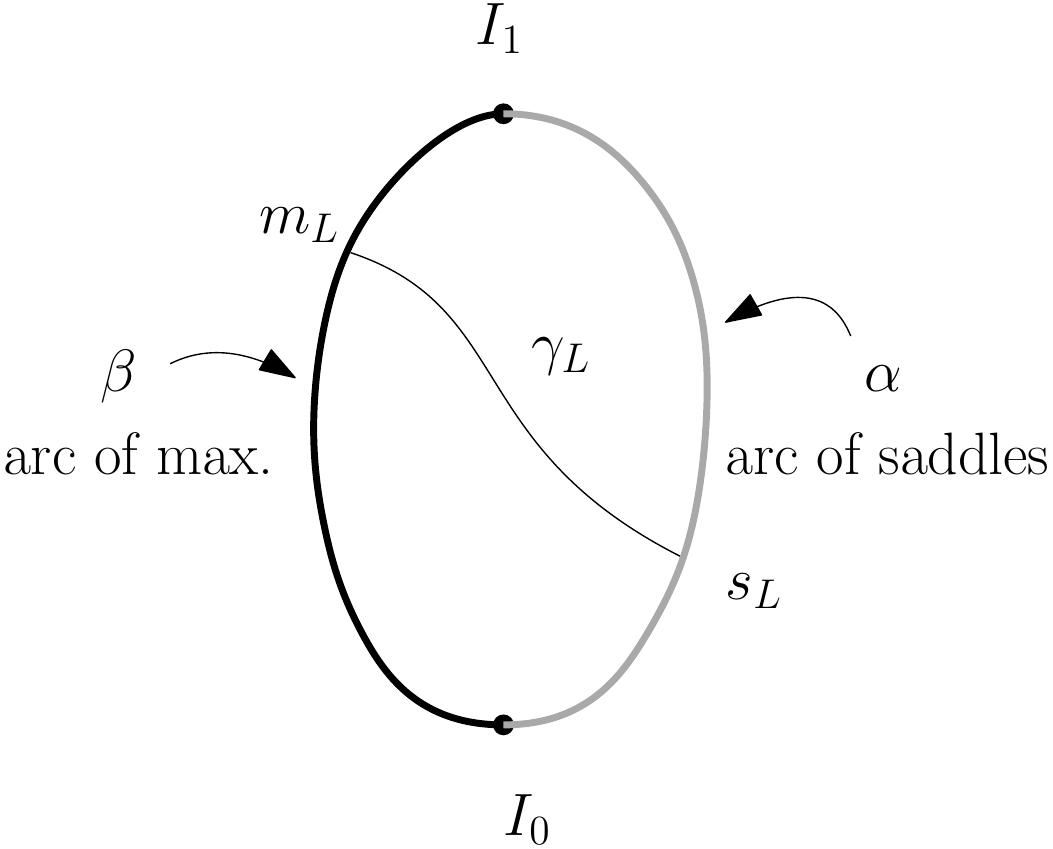}
\caption{ A simple loop max-saddle.}\label{loop-nb}
\end{figure}
 \end{center}

\proof We may only consider the case where $\beta$ is an arc of maxima. Let 
$\ep$ be a small positive number. 
For every $m_L\in \beta$,  consider the stable manifold 
$W^s(m_L, \nabla_\L f)$, truncated at the level $f=f(s_L)-\ep$.  
 Then slightly
enlarge this surface to a disc $ D_{m_L}\subset L$ in order to have a smooth boundary and to contain
 an $\ep$-neighbourhood of the truncated stable manifold of $s_L$. If $\ep$ is small enough, $D_{m_L}$
avoids all other connected components of the contact locus. 
If this construction is performed smoothly with respect to $m_L$ with a fixed $\ep$,
 the desired $N$ is $\cup_{m_L\in\beta}D_{m_L}$.

There exists a proper disc $S\subset N$, transverse both to the leaves of $\ell$ and the level sets of $f$,
that contains all segments $\ga_L$ (item (3) in Definition \ref{def-simple}.) This $S$ can be endowed with coordinates
$(x,z)$ where $f(x,z)= z$ and the pair $(x,z)$ allows us to use the normal form of inflections. Finally, one finds 
the $y$-coordinate by slicing $N$ transversely to $S$ so that the maximum 
of $f_{\vert \{x=x_0, \ell(x_0,y,z)=\ell_0\}}$  is quadratic non-degenerate and located in $\{y=0\}$.
In Morse theory, such extension of coordinates is named a \emph{suspension} \cite{h-w}.

We are able to  solve the problem of isotopy for a pair of functions like $(f_{\vert S}, \ell_{\vert S})$ since the two-dimensional problem is known  by the contractibility of  Diff$(D^2 \text{ rel. } S^1)$ \cite{smale}. An isotopy of $S$
extends to $N$, relatively to a small neighbourhood of $\p N$ by suspension. 
Nothing is changed 
in $(S^2\times[0,1])\smallsetminus N$. 
The proof is now complete.\footnote{ This proof is inspired by my proof 
of a theorem of Morse about the cancellation of a pair of critical points of one real function \cite{l-cras}.}\bull

On the opposite, creating a simple loop requires no condition. 
It can be realized in any open set of the ambient manifold where the
two functions $f$ and $\ell$ have no contact points. The details are left to the reader.

The next elementary isotopy will follow the same idea of suspension. This deals with the configuration of contacts
named {\it saddle-center-saddle}. Center stands for minimum or maximum; here, we only consider the case of a 
minimum. 
The case where the center is a maximum is similar.\\

\begin{defn} \label{scs-conf}
The saddle-center-saddle configuration is the following. The three contact points are in the same leaf $L$
of $\L$; the center is a minimum $m_L$ whose cone $C_\L(m_L)$ contains two saddles $s_L$ and $s'_L$
in its boundary; there is exactly one $\L$-gradient line from $m_L$ to $s_L$, and similarly from $m_L$ to $s'_L$.
So, $f(s_L)=f(s'_L)$. 

Let $\al, \beta, \al'$  denote the contact arcs containing $s_L, m_L, s'_{L}$ respectively. Here are the 
main requirements: 
\begin{enumerate}[\rm (1)]
\item Let $L'$ be a leaf close to $L$ crossing $\beta$ at $m_{L'}$ with $f(m_{L'})<f(m_L)$. Then $f(s_{L'})<f(s'_{L'})$.
\item Let $L''$  be a leaf close to $L$ crossing $\beta$ at $m_{L''}$ with $f(m_{L''})>f(m_{L})$.
Then $f(s_{L''})>f(s'_{L''})$.
\item The interior of the ball $B_\L(m_L)$ meets no other contact arc than an arc of $\beta$.\footnote{See subsection
\ref{notation-ball} for notation $B_\L(m_L)$.}
\end{enumerate}
Here, it is meant that $\{s_{L'},s_{L"}\}\subset \al$ and $\{s'_{L'},s'_{L"}\}\subset \al'$.
\end{defn}

\begin{center}
 \begin{figure} [h]
\includegraphics[scale =.7]{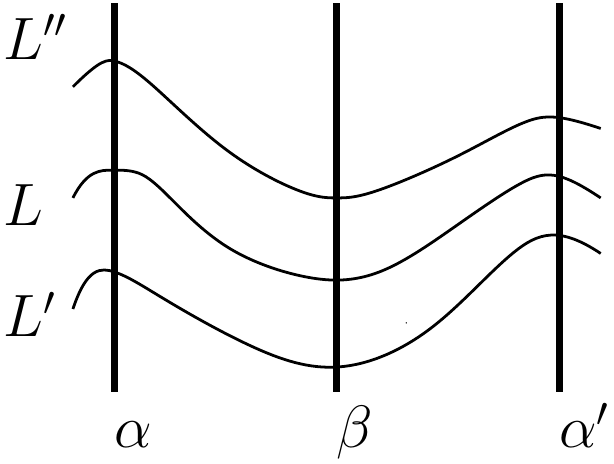}
\caption{ The leaves are denoted on the left. Vertically are the three contact arcs.}\label{planar-fig-nb}
\end{figure}
 \end{center}
 
\begin{center}
 \begin{figure}[h]
\includegraphics[scale =.7]{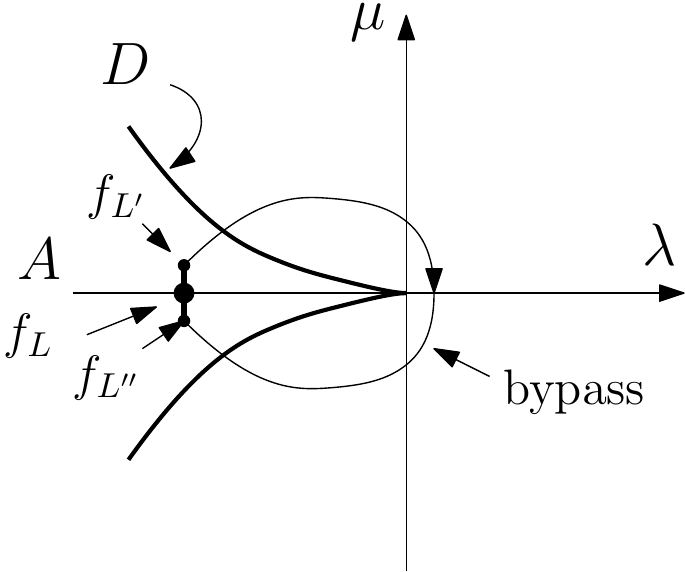}
\caption{Discriminant locus of the singular function $x^4$ or $x^4-y^2$.}\label{cusp-fig-nb}
\end{figure}
 \end{center}

 Note that, up permuting $\al$ and $\al'$, (1) and (2) above-mentioned 
 follow from the excellence of the pair $(f,\ell)$; moreover,
it is easily seen by (3) that $s_L, \ m_L,\text{ and } s'_L$ have the same sign as contact point.
The planar figure of this configuration is represented in Figure \ref{planar-fig-nb}. 
One recognizes---up to sign---the 
singularity $x^4$, named {\it cusp},
whose  versal unfolding---up to conjugation source{\tiny x}target---reads
\begin{equation}\label{fold1}
(x, \la,\mu)\mapsto -(x^4+\la x^2+\mu x)
\end{equation}
where $\la$ and $\mu$ are two real parameters. As $m_L$ is a minimum, the suspension  consists of adding a {\it positive square} in the added variable $y$. Hence, it reads 
\begin{equation}\label{fold2}
(x,y, \la,\mu)\mapsto -(x^4+ \lambda x^2+\mu x) + y^2.
\end{equation}
The critical points of $f$ restricted to a leaf close to $L$ are given by $y=0, 4x^3+2\la x+\mu=0$.
 Up to  positive coefficients, the equation of this discriminant locus reads 
\begin{equation}\label{discriminant}
4\la^3+ 27 \mu^2=0.
\end{equation}
The $(\la,\mu)$-space, denoted by $\R_{(\la,\mu)}$,
is stratified,
 apart from the origin, by the discriminant locus $D$ of the \emph{fold},  and the half axis $A=\{\la<0\}$, the set of 
parameters for which the corresponding function has two
 critical points with the same value. This is exactly the case of $f_{L}$.
 The deformation from $L'$ to $L''$ corresponds to a path in $\R_{(\la,\mu)}$ crossing $A$ transversely 
 (Figure \ref{cusp-fig-nb}.)

In the setting saddle-center-saddle, we have the 3-ball $B_{\L}(m_L)$ that we can extend up to the level 
$\{f= f(s_L)+\ep\}$. The closure of this extension contains a piece of the unstable manifolds
$W^u(s_L, \nabla_\L f)$ and $W^u(s'_L, \nabla_\L f)$. Then, we take a neighbourhood $N$ of this closure.
And hence, there are coordinates $(x,y,z)$ on $N$ such that the surface $S:=\{y=0\}$ is the locus of the minimum of $f$ 
on the curve 
in $N$ defined by $(x,\ell)=(x_0,\ell_0)$ where the pair $(x_0,\ell_0)$ ranges in a convenient 2-dimensional domain.

\begin{center}
 \begin{figure}[h]
\includegraphics[scale =.7]{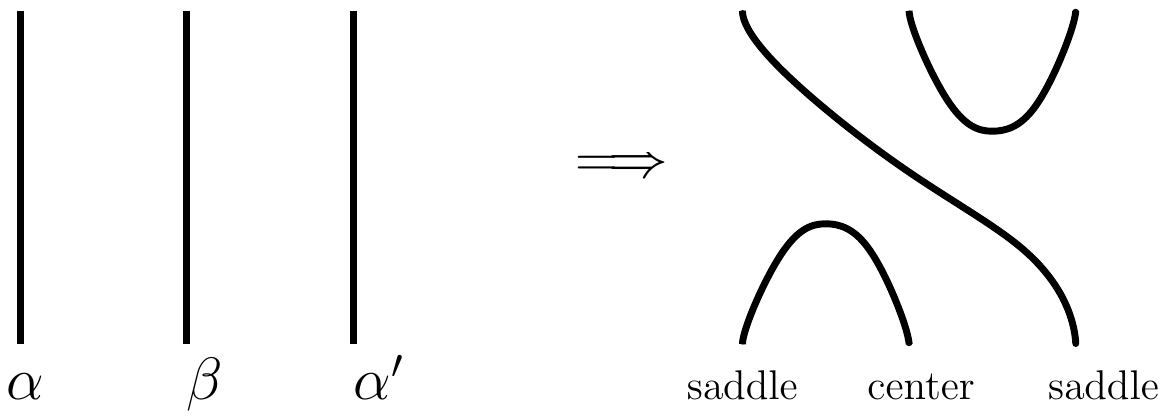}
\caption{Bypass in configuration saddle-center-saddle.}\label{bypass-fig-nb}
\end{figure}
\end{center}

\begin{defn}
A \emph{bypass} of the cusp singularity is any path of functions as in {\rm Figure \ref{cusp-fig-nb}.}
\end{defn}
Figure \ref{bypass-fig-nb} represents the modification of the contact locus along such a path.

\begin{prop} \label{scs-prop}
In the above-mentioned setting, there is an isotopy supported in $N$ that, at time $1$,
carries $\ell$ to a function $\ell'$ whose one-parameter family of leaves realizes a bypass of the cusp singularity.
\end{prop}

After these first two elementary isotopies, we present more shortly one useful variation of the latter.

\begin{center}
 \begin{figure}[h]
\includegraphics[scale =.7]{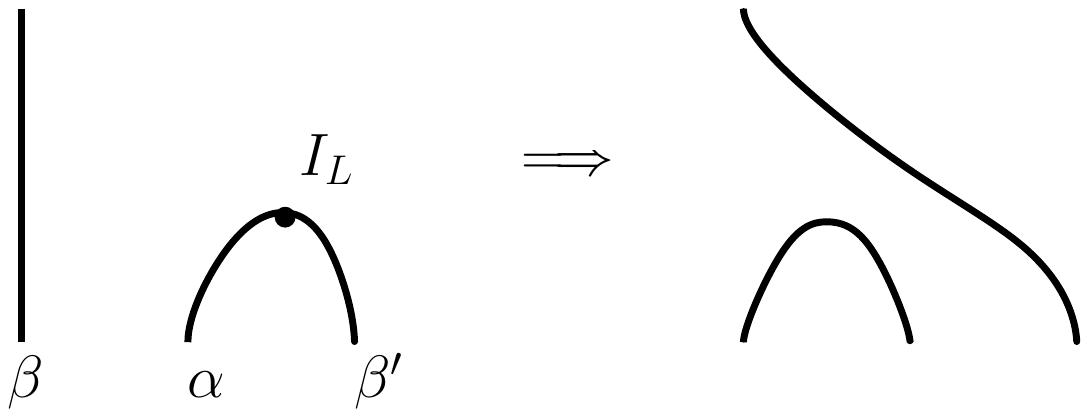}
\caption{\small Min-saddle-min configuration.}\label{bypass2-nb}
\end{figure}
\end{center}

\begin{defn}\label{msm-conf}
The configuration \emph{min-saddle-min} is the following. A quadratic contact $m_L$ of index $0$, located in the leaf 
$L$, has a cone $C_\L(m_L)$ whose boundary has an inflection 
$I_L$  of the same sign as $m_L$.
Let $\beta$ be the contact arc containing $m_L$. It is assumed:
\begin{enumerate}[{\rm 1)}]
\item The ball $B_{\L}(m_L)$ contains no other contact arcs but a sub-arc of $\beta$.
\item $I_L$ is the cancellation point of a pair saddle-min, that is the common upper bound of an arc $\al$ of saddles 
and an arc $\beta'\neq\beta$ of minima (Figure \ref{bypass2-nb}, left-hand side.)
\end{enumerate}
\end{defn}

Therefore, there is an $\L$-gradient line from $m_L$ to $I_L$. So,
 if $m_{L'}$ is just below $m_L$ on $\beta$, 
 then there is a saddle $s_{L'}$ connected to $m_{L'}$ by an $\L$-gradient line, 
 and hence, the other branch of the stable manifold $W^s(s_{L'}, \nabla_\L f)$  comes from a minimum  $m'_{L'}$
 distinct from $m_{L'}$. 
 
 If one forgets the so-called suspension by $+y^2$, it is clear that the new configuration deals with 
 the cusp singularity $x^4$. Instead of crossing the equality of  two critical values, one crosses the cancellation
 of a saddle with a center.
 
 \begin{center}
 \begin{figure}[h]
\includegraphics[scale =.7]{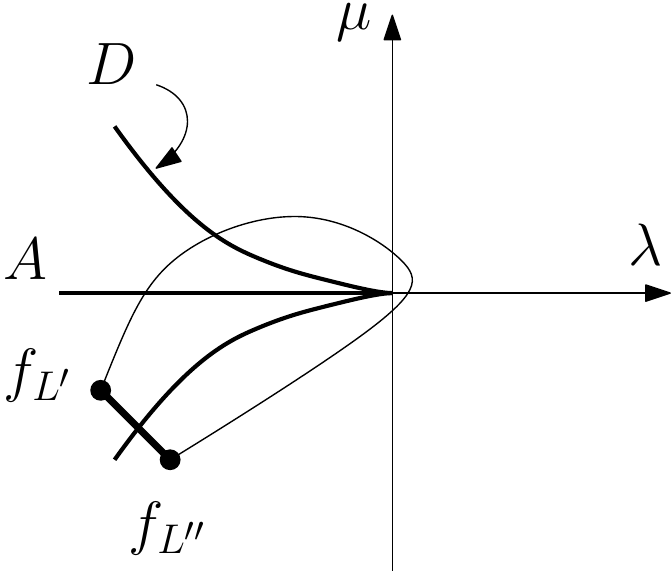}
\caption{\small Another bypassing of the cusp singularity.}\label{exchange-nb}
\end{figure}
\end{center}
In this setting, bypassing the cusp singularity has the following effect on the contact curves, as shown in Figure 
\ref{bypass2-nb}.

\begin{prop} \label{msm-prop}
In the configuration min-saddle-min, there is an isotopy supported in a neighbourhood $N$
of the ball $B_\L(m_L)$ that, at time 1, carries $\ell$ to a function whose contact curves 
in $N$ have the behaviour presented in the right hand side of {\rm Figure \ref{bypass2-nb}.}
\end{prop}

We now turn to \emph{isotopies along saturated sets}.  
\begin{defn}  \label{satur}
Let $x\in S^2\times [0,1]$ 
and $z_0>f(x)$.
The  \emph{ascending saturated set} of $x$ up to level $z_0$,
denoted by $S_{z_0}(x)$, is the minimal closed subset of $\{f\leq z_0\}$ that contains $x$ and satisfies:
\begin{enumerate} [\rm (i)]
\item For every $y\in S_{z_0}(x)$, the positive orbit of $y$ under $\nabla_\L f$, truncated at level $z_0$, is included in
$S_{z_0}(x)$.
\item Let $C$ denote the contact locus of $(\F,\L)$.  If $y\in (C\cap S_{z_0}(x))$,  every arc of $C$ ascending 
 from $y$ and truncated at level $z_0$ 
 is contained in $S_{z_0}(x)$.  
\end{enumerate}
\end{defn}
Of course, there is an analogous definition of a \emph{descending satured set}.
Such a subset is stratified in an obvious way. By excellence assumption, $S_{z_0}(x)$ 
contains finitely many \emph{accidents}, 
meaning  level sets where the stratified type changes. 
From one accident to the next one, one proves
that $S_{z_0}(x)$ is surrounded  by a fundamental system of \emph{collapsible domains} that, in our setting, 
are defined as follows.

\begin{defn}\label{collapse} 
A \emph{collapsible domain} $K$ is a 3-ball  whose boundary $\p K$ is the angular union of two discs, the \emph{lower} 
boundary $\p_{lo} K$ and the \emph{horizontal} boundary
 $\p_h K$, fulfilling the next conditions:
 \begin{enumerate}[\rm (i)]
 \item $\p_h K$ is contained in $\{f= z_0\}$.
 \item $\p_{lo}K$ \rm (resp. its interior) is contained in $\{f\leq z_0\}$ (resp. $\{f<z_0\}$).
 \item There exists a field of directions on $K$ transverse to the $f$-level sets
 and also transverse to $\p_{lo}K$.
 \end{enumerate}
\end{defn}

\begin{lemme} \label{collapse-lem} Let $K$ be a collapsible domain  which
contains  $S_{z_0}(x)$ and whose lower boundary  is disjoint from $S_{z_0}(x)$. For every $\ep>0$,
there exists an isotopy $(\Phi^t)_{t\in[0,1]}$, supported in the interior of $K$, such that:
\begin{enumerate}[\rm(1)]
\item $\Phi^0= Id$ and $\Phi^1(S_{z_0}(x))$ is contained in $\{z_0-\ep<f\leq z_0\}$;
\item for every $t\in[0,1]$ and every contact point $y$ of the pair $(f,\ell)$, its image
 $\Phi^t(y)$ is a contact point of the pair $(f,\ell\circ(\Phi^t)^{-1})$,
with the same sign and Morse index;
\item if $\ga$ is an $\L$-gradient arc, then $\Phi^t(\ga)$ 
is transverse to the $f$-level sets.\footnote{ See more details in
\cite[chap. III §4]{lauden}.}
\end{enumerate}
\end{lemme}

Such an isotopy is said to be an \emph{isotopy along a saturated set}.\\
  
\begin{rien}{\sc Application.}\label{saddle--text}
{\rm Here we  give an application of Proposition \ref{scs-prop}, namely, the cancellation of a pair of saddles, 
that are of the type $X$.  We specify the definition of an $X$-saddle when the pair $(f,\ell)$ is excellent.

\begin{defnt} 
A saddle $s$ in a contact arc $\al$ is said to be of type $X$, or an \emph{$X$-saddle},
if $s$ is the upper bound on $\al$ of an interval of saddles of type $\la$ (resp. $Y$) and the lower bound of an interval of 
saddles of type $Y$ (resp. $\la$).
\end{defnt}

The following facts are easily checkable when the pair $(\F,\L)$ is excellent and the leaves are 2-spheres.

\begin{lemme} \label{X-pair}

{\rm (1)} If $s$ is an $X$-saddle there is another $X$-saddle $s'$ in some contact arc $\al'$ 
interacting with $s$, meaning that  the two saddles
are located in the same connected component $\hat C_0$ of the intersection $ L_0\cap F_0$ of a leaf 
$L_0:=\{\ell= u_0\}$ and a level set $F_0:=\{f=z_0\}$.

{\rm (2)} Both saddles have the same sign and $\hat C_0$ looks like the bold line on 
Figure \ref{saddle-nb}. 

{\rm (3)} $\hat C_0$ is made of two closed singular curves $C_0$ and $C'_0$ in $F_0$ 
 meeting in $s$ and $s'$ only.\footnote{There are two possibilities for the pair $(C_0,C'_0)$: either it is the limit 
 of the dotted curve from Figure \ref{saddle-nb} or it is the limit of the dashed curve.}
\end{lemme}

\begin{center}
   \begin{figure}[h]
\includegraphics[scale =.4]{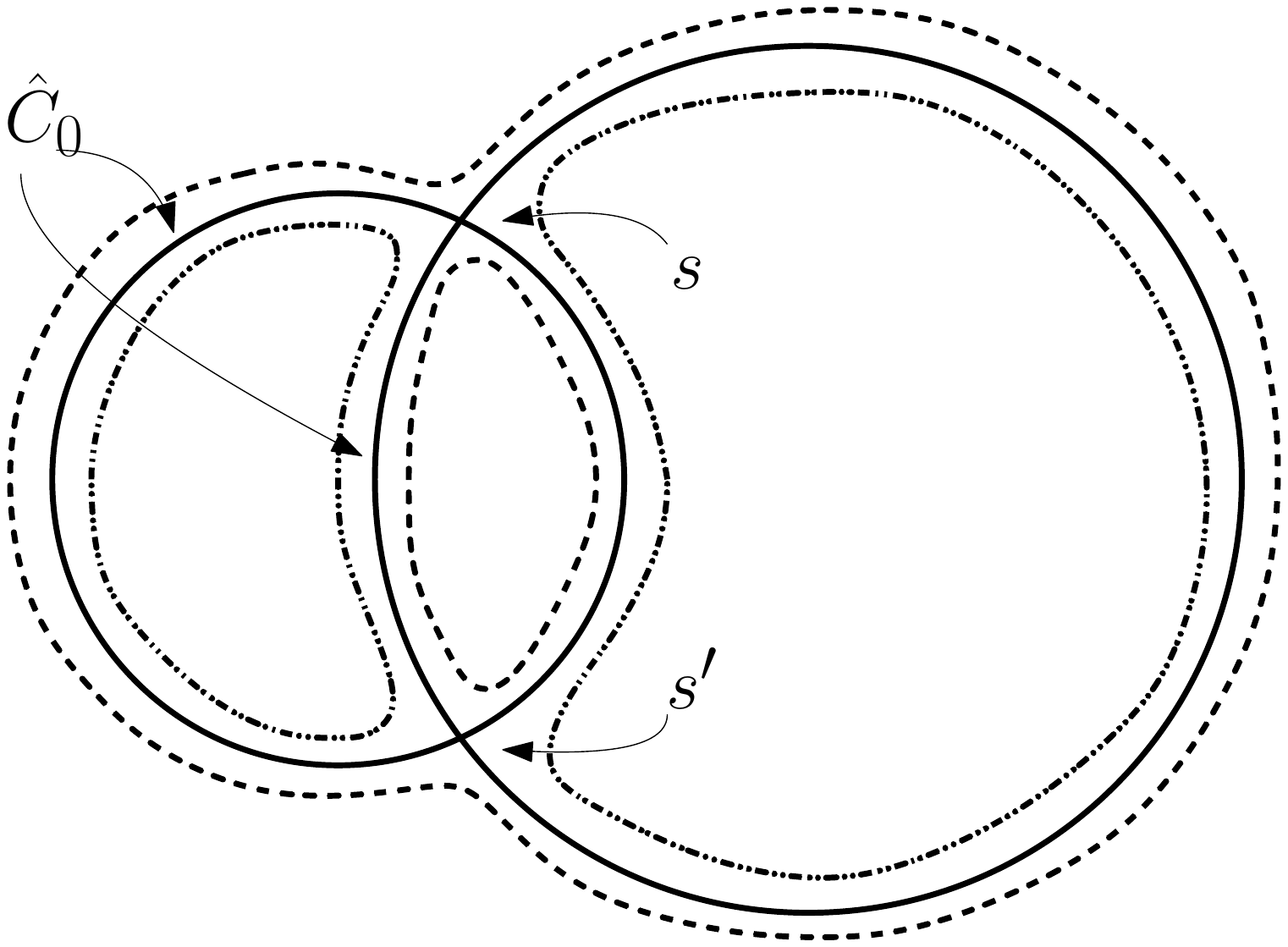}
\caption{\small $\hat C_0$ denotes the connected component of $s$ in the level curve $F_0\cap L_0$. The dashed (resp. dotted) curves lie in $L_0$ and in level sets above (resp. below) $z_0$. }\label{saddle-nb}
\end{figure}
 \end{center}

\proof
It almost clear that a second saddle in $\hat C_0$ is necessary. Indeed, if not, the germ of $\L\cap F_0$
$\hat C_0$ is stable
in the sense that, up to isotopy, it does not depend on small variations of $z_0$. Therefore, $s$ could not be of type $X$.

\begin{center}
   \begin{figure}[h]
\includegraphics[scale =.6]{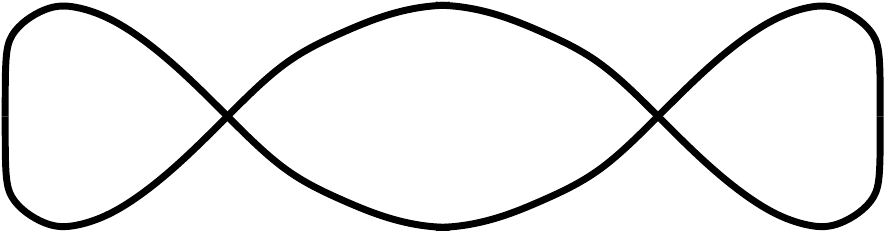}
\caption{\small The unique example of a chain in $S^2$. Uniqueness up to isotopy does not hold in $\R^2$.}\label{saddle-2-nb}
\end{figure}
 \end{center}

By excellence, $\hat C_0$ cannot have three saddles (or double points.) Then, it contains exactly two saddles. Up to isotopy of $F_0$, 
there are very few configurations for $\hat C_0$: 
\begin{enumerate}[(a)]
\item The one represented in Figure \ref{saddle-nb}. 
\item A second one that is the image---named \emph{a chain}---of a smooth closed cuve through an immersion with two double points (Figure \ref{saddle-2-nb}.)
\end{enumerate}

The rule for the co-oriention of $\hat C_0$ in $F_0=S^2$ near the double points leaves only these two cases. In the 
same way, one checks that in case (a) the two contact points have the same sign. 

In case (b), whatever the signs of contacts the saddles cannot be $X$-saddles. \bull

 \begin{center}
 \begin{figure}[h]
\includegraphics[scale =.6]{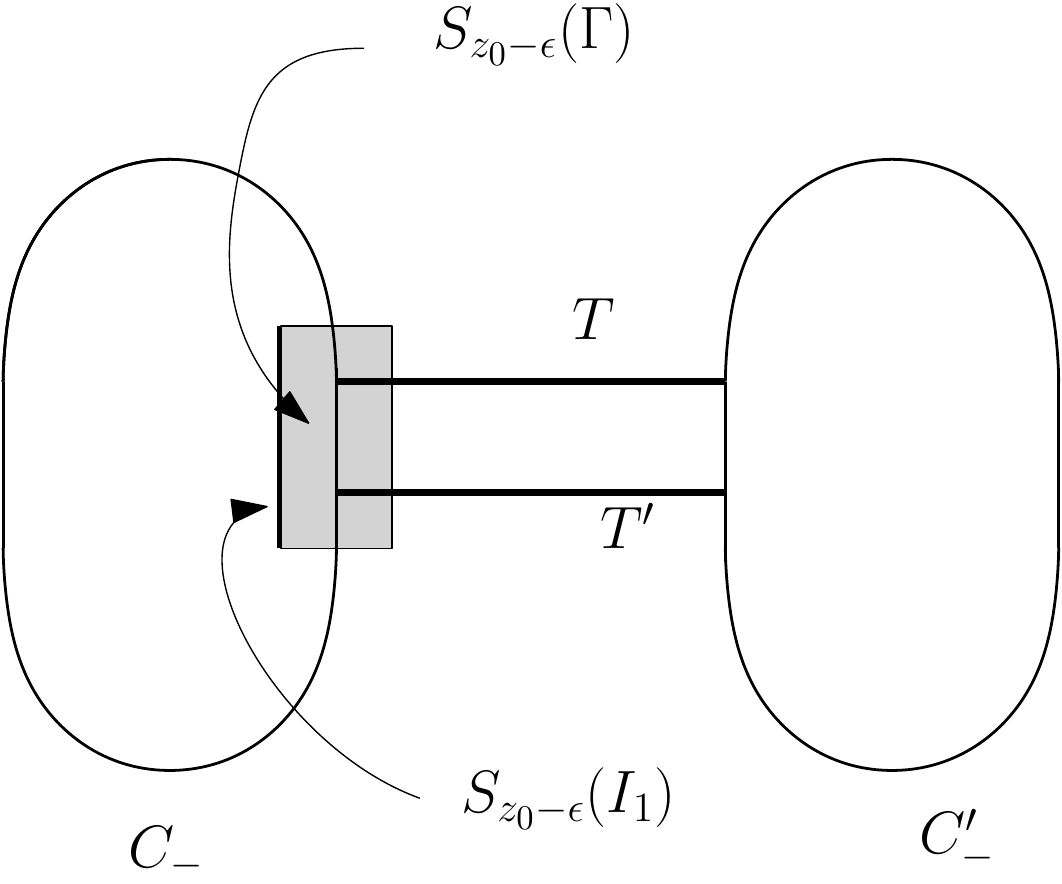}
\caption{\small The gray rectangle represents the domain marked by the ascending saturated set of $\Ga$ in the level set $z_0-\ep$ .}\label{X-nb}
\end{figure}
\end{center}

Denote by $\al'$ the arc of saddles passing through $s'$.  Let $s(t)$ and $s'(t)$ be two local 
  parametrizations of $\al$ and $\al'$, with $s(0)=s, \ s'(0)=s'$, such that $\ell(s(t))=\ell(s'(t))$ for every $t$ close to $0$. Say $f(s(t))$ and
  $f(s'(t))$ are increasing with $t$. As the pair $(f,\ell)$
  is excellent, the two $f$-velocities at $t=0$ are different. We continue this analysis by specifying that the two contacts 
  are positive. 

  For $\ep$ small enough, let $C_-$ and $C'_-$  be the two smooth closed curves in the level set 
  $\{f= z_0-\ep\}$ corresponding to $C_0$ and $C'_0$ respectively by the $\L$-gradient flow lines on $L_0$.
  Let $T$ and $T'$ be the lower boundaries of the  saturated sets descending 
  from $s(0)$ and $s'(0)$ respectively (Figure \ref{X-nb}).
}
\end{rien}

\begin{prop}\label{cancel-X}
There is an isotopy supported in $\{z_0-2\ep \leq f\leq z_0+\ep\}$ whose effect on $\ell$ is
to decrease by one the number of pairs of contact points with $\F$ that are of type $X$.
\end{prop}
\proof 
The idea is to create a configuration saddle-center-saddle which contains $s(0)$ and $s'(0)$. Since the center is 
missing, one starts by introducing a simple loop $\Ga$ of type saddle-min
 whose inflection points $I_0$ and $I_1$ satisfy 
$z_0-2\ep<f(I_0)<f(I_1)<z_0-\ep$ and so that its upper boundary in $\{f=z_0-\ep\}$ is a position 
as shown in Figure \ref{X-nb}.

\begin{center}
 \begin{figure}
 \includegraphics[scale =.7]{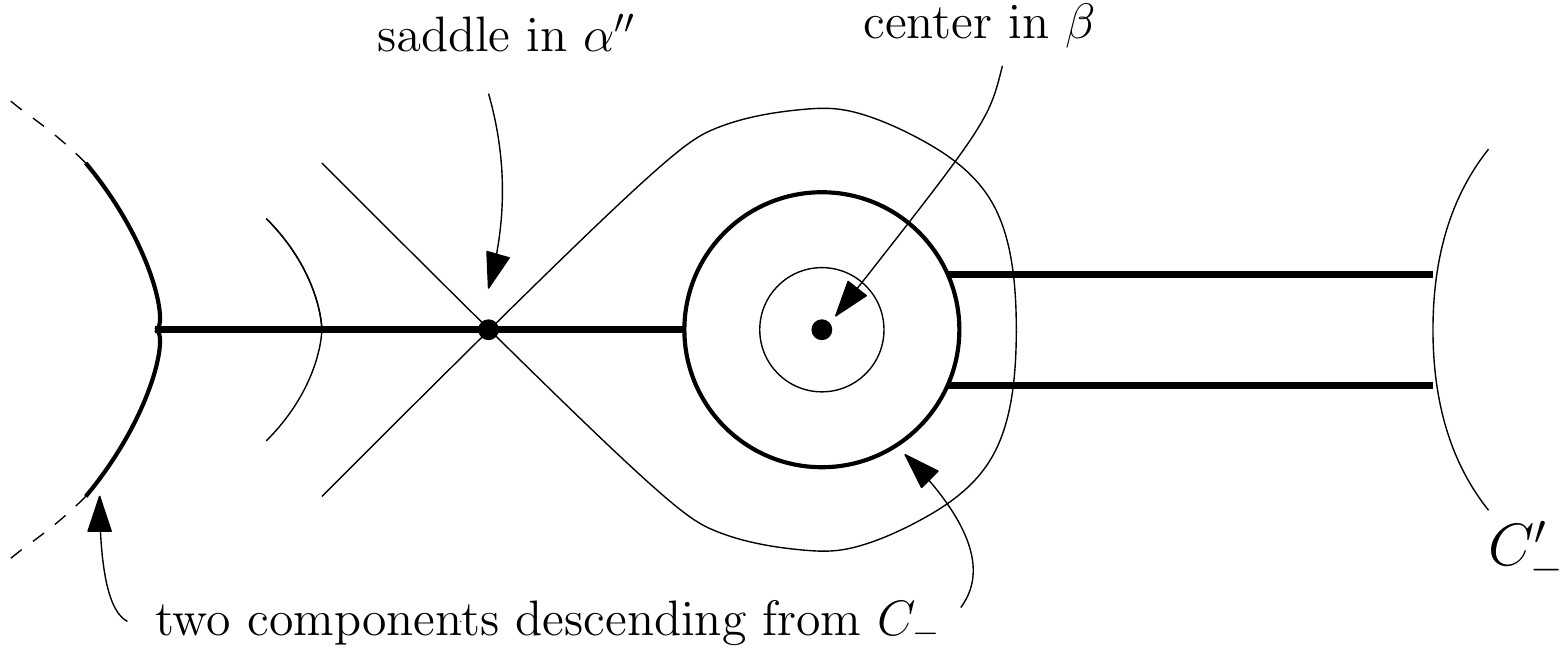}
\caption{\small The very thick lines are lower boundaries of descending saturated sets of saddles $s(0)$, $s'(0)$
and $\al''\cap \{f=z_0-\ep\}$.}\label{X-min-nb}
\end{figure}
\end{center}

The next operation is an isotopy along the satured set $S_{z_0+\ep}(I_1)$. By application of Lemma \ref{collapse-lem},
one can collapse this domain above $z_0$. After this isotopy, still denoting by 
$\ell$ and $\Ga$ the carried function
and simple loop, one sees a configuration saddle-center-saddle as shown in Figure \ref{X-min-nb}. Then Proposition
\ref{scs-prop} is available that makes the desired cancellation of a pair of $X$-saddles.

There is a new arc $\al''$ of positive saddles which crosses $\{f=z_0-\ep\}$---this is the arc of saddles in the simple loop 
$\Ga$ after the just above mentioned isotopy; the arc of minima is denoted by $\beta$. 
One checks that $\al''$ 
contains only
saddles of type $\la$. Now we can apply Proposition \ref{scs-prop}. This cancels the pair $(s,s')$ of $X$-saddles 
without creating a new one and  Proposition \ref{cancel-X}
 is proved.
 
 \bull
 
 The outcome of Proposition \ref{cancel-X} is the cancellation of every $X$-saddle.

 \section{Dehn modification} \label{dehn}
 
After Thom's transversality arguments that generate generic properties,  the Dehn modification, also called
Dehn surgery,  is the main tool to prove Theorem (A), and already regarding its 1-connected version Theorem (A'). 
We are going to speak of this notion in our setting only: the ambient manifold is $M:=S^2\times [0,1]$,
provided with two foliations, respectively as the level sets of two functions, $f$ and $\ell$, 
that are constant on each component 
 boundary components and 
whose critical sets are empty.

As an application, we prove that a suitable series of Dehn modifications allows us to 
move the $Y$-saddles 
above the $\la$-saddles. 

\begin{defn} {(\sc Dehn twist and Dehn modification.)} \label{dehn-def}${}$

\nd{\rm 1)} A \emph{left Dehn twist} 
 on the standard annulus $\mathbb A:= S^1\times [1,2]$ is a diffeomorphism $\tau$
that, in polar coordinates, reads 
\begin{equation}\label{twist-formula}
(\theta, r)\mathop{\longmapsto}^\tau(\theta +\vp(r-1),r)
\end{equation}
where $\vp$ is a smooth non-decreasing function from the value $\vp([0,\ep])= 0$ to $\vp([1-\ep, 1])= 2\pi$.

The \emph{right} Dehn twist is the inverse $\tau^{-1}$ of the left one.
\smallskip

\nd {\rm 2)} Let $A$ be an annulus embedded in $M$ and parametrized by $\mathbb A$; for $i= 1,2$, $\p_iA$ is parametrized by $\p_i\A:= S^1\times\{i\}$. 
Let $\check M$ be the manifold obtained by cutting $M$ along the interior of $A$; it is provided with a singular boundary 
made of two lips $A^+$ and $A^-$ corresponding to the coorientation of $A$. The \emph{Dehn modification} of $M$
along $A$ is the manifold $M_\tau$ obtained from $ \check M$ by gluing $x\in A^+$ to $\tau(x)\in A^-$.
\end{defn}

 \begin{center}
 \begin{figure}[h]
 \includegraphics[scale =.6]{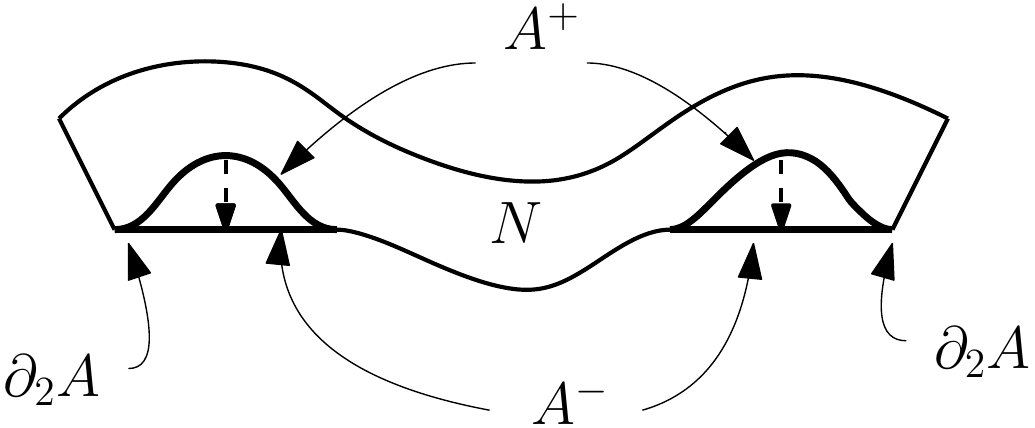}
\caption{\small The two dotted vertical arrows indicate the gluing 
that makes $M_\tau$.}\label{dehn-nb}
\end{figure}
\end{center}

 In our setting, $A$ is an annulus, embedded 
 in a level set $F_0$ of $f$ transversely to $\L$; the leaves of $\L\cap A$ are closed curves isotopic to each 
 component of $\p A$.
 The identification with $\mathbb A$ is chosen so that these curves are parametrized by the circles of $\A$.
 In this setting we shall say that the Dehn modification is \emph{adapted} to the pair $(f,\ell)$ or $(\F,\L)$.
 By construction, the functions $f$ and $\ell$ are carried to some uniquely defined  functions $f_\tau$ and $\ell_\tau$. \\
 
 \begin{important-rem}\label{important}   In the previous setting, a Dehn modification along $A$ keeps the contact locus of the pair
 $(f,\ell)$ invariant. {\rm Indeed, the support of this modification is away from the contact locus.}
 \end{important-rem}

 \begin{defn}  \label{disc-space} Let $A$ be an annulus, parametrized by $\A$, in a level set $F_0=\{f=z_0\}$.
 Let $p_A\in S^2$ such that the vertical arc $\{p_A\}\times[0,1]\subset S^2\times[0,1]$ avoids $A$ 
 and intersects the unique disc in $F_0$ that contains  $A$ and is bounded by $\p_2A$.
 
 Let $\mathcal D_A$ denote the space of smooth oriented 2-discs embedded in $M$  that contain $A$
as oriented annulus, 
have algebraic intersection $+1$ with $\{p_A\}\times[0,1]$
 and have $\p_2A$ as a boundary. 
  \end{defn}

 \begin{prop}\label{diffeo} With this notation, 
   the following holds true.
 \begin{enumerate}[{\rm 1)}]
  \item Every $D\in \mathcal D_A$  defines a unique  
  diffeomorphism $\psi_D: M_\tau \to M$, up to isotopy. 
 \item For every pair $(D_0,D_1)$ of elements in $\mathcal D_A$, then $ \psi_{D_1}$ is isotopic to $ \psi_{D_0}$.

 \item There are two discs  $D_0$ and $D_1$, elements of $\mathcal D_A$, such that 
 $(\psi_{D_0})^* df=df_\tau, \text{ and } (\psi_{D_1})^*d\ell =d\ell_\tau$.

\item The functions $f$ and $\ell$ are isotopic if and only if the functions $f_\tau$ and $\ell_\tau$ are isotopic.
 \end{enumerate}
  \end{prop}
  
  \proof 1) Since $A$ is an imposed collar of $\p D$ for every  $D\in \mathcal D_A$, it is natural to say that, 
  abstractly, $D$ is the disc of radius 2
  in the Euclidean plane.
  This disc may be considered as embedded in $M$
  or in $M_\tau$ as well. Firstly, we consider $D$ in $M_\tau$. One endows $D$ with a 3-dimensional collar 
  $N\cong D\times[0,1]$
  on the side of $A^+$, that is $A^+\subset D\times\{0\}$. So, $N$ already existed  in $\check M$ and 
   $M$.

  One defines a diffeomorphism $\rho: D\to D$ by $\rho(x)=\tau(x)$ if $\Vert x\Vert\in [1,2]$ and 
  $\rho(x)=x$ if $\Vert x\Vert\leq 1$. One chooses an isotopy $(h_t)_{t\in [0,1]}$ among the diffeomorphisms
  of $D$ from $\rho$  for $t= 0$
  to the identity of $D$ for $t= 1$, with an extra requirement:
  \begin{equation}\label{norm}
  \Vert h_t(x)\Vert =\Vert x\Vert\text{ for every }t\in[0,1].
  \end{equation} 
  Then, on defines $\psi_D:M_\tau \to M$ by 
  \begin{equation} \label{formule}
  \left\{
  \begin{array}{ll}
  \psi_D(y)&=y \text{ if } y\in M\smallsetminus N \\
  \psi_D(x,t)&= (h_t(x),t)  \text{ if } (x,t)\in N.
  \end{array} \right.
  \end{equation}
  Since all choices are made in contractible spaces, $\psi_D$ is uniquely defined up to isotopy.
  \begin{center}
 \begin{figure}[h]
\includegraphics[scale =.6]{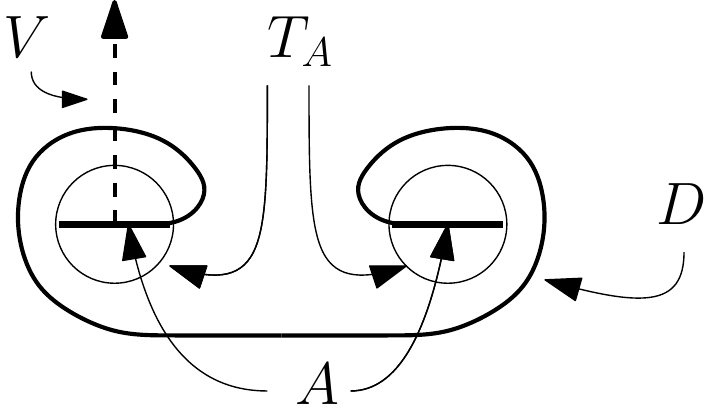}
\caption{\small  Example of meridional twist. Here, $w_A(D)= -1$.
}\label{winding-nb}
\end{figure}
\end{center}
  
   2) Usual tools of 3-dimensional topology yield the following topological result.\footnote{ Unfortunately, I have no 
   reference for this exercise. One has to use the classical method of the \emph{innermost} intersection curve for two 
   discs  in $\mathcal D_A$ with the same winding number $w_A$. }\\
   
   \parbox{15cm}
   {\it Each 2-disc  $D\in \D_A$ has a \emph{winding number} around $A$, $w_A(D)\in \Z$, that is the algebraic 
   intersection of $D$ with a \emph{vertical} arc $V$ in $S^2\times [0,1]$ starting from the interior of $A$ and ending in $S^2\times \{1\}$.
   Moreover,  $w_A$ provides a bijection from $\pi_0(\D_A) $ onto $\Z$.} \\
   
   \rm Here, the linking condition imposed to 
   $D$ for being in $\mathcal D_A$ plays an important role.
  This statement also holds in $M_\tau$. Moreover,
   one can think of $A$ as a flat annulus in a solid torus $T_A\cong S^1\times D^2$, each ray of $A$ being a diameter of 
   $\{\theta\}\times D^2$.
    An exterior collar of $\p T_A$ has coordinates $(\theta, \theta', t)$ with $(\theta',t)\in \p D^2\times [0,1]$. 
    This allows one to define a \emph{meridional Dehn twist}
    in the coordinates $(\theta',t)$ independant of $\theta$.\footnote{ Use formula (\ref{twist-formula}) replacing 
     $(\theta,r)$  
    with $(\theta',t)$.} Let us denote it by $\mu$.
    For $D\in \D_A$, one has 
    \begin{equation} \label{mu-formula}
    w_{A}(\mu(D))= w_{A}(D)+1.
    \end{equation} 
    This holds true on $M$ and $M_\tau$  as well because the support of $\tau$ and $\mu$ are disjoint.
    
    Using that translations on the torus are commuting, one has for every $D\in \D_A$:
    \begin{equation}
    \psi_{\mu(D)}= \mu\circ \psi_D\circ\mu^{-1}
    \end{equation}
    Moreover, $\mu $ is isotopic to the identity if $A$ is let free to rotate in $T_A$ around its core.
    rotating $A$ does not change $M_\tau$. Therefore, $\psi_{\mu(D)}$ is isotopic 
    to $\psi_D$ which proves what was desired.\\

  3) The disc $D_0$ is the unique disc in the level set $F_0$ fulfilling the requirement $\p D_0= \p_2A$.
  One can choose its collar $N\cong D_0\times[0,1]$ to be foliated by level sets of $f$. Then  formula (\ref{formule})
  shows that
  $(\psi_{D_0})^*f$ is sensitive neither to $\tau$ nor to the isotopy $(h_t)$. So, $(\psi_{D_0})^*f=f$ is clear.
  
  For $D_1$, this is slightly more subtle. There is a disc $\De$ which is  bounded by $\p_1A$ in some $\L$-leaf 
  and intersects $\{p_A\}\times[0,1]$ positively.
  The union $A\cup \De$ is angular along $\p_1A$. But its smoothing is easy:
   Take a collar $N_\De$ of $\De$ in the  direction outward to $A$ along $\p_1A$ and foliated by $\L$.
   Then,
    cap  $\p_1A$ in $N_\De$ with a disc $\tilde \De$ tangent to $A$ along $\p_1A$ and so that 
    $\L$ induces a foliation by circles with one center. Now, we take $D_1= A\cup \tilde \De$. 
    
    The last step consists of choosing a collar $N_1= D_1\times [0,1]$ of $D_1$ in $M$, that 
     exists in $\check M$ (on the side of $A^+$) so that $\L$ induces 
     on each $D_1\times \{t\}$ a foliation by circles. Therefore, 
    by formula (\ref{formule}) the isotopy  is compatible with $\L$. And hence, $\psi_{D_1}(\L_\tau)=\L$.\\
   
    4) This is a formal consequence of items 2) and 3). \bull
  
  We have seen that a Dehn modification adapted to $(\F,\L)$ allows us to carry 
  some structures, like $f$ or $\ell$, from $M$ to $M_\tau$. But an isotopy, for instance along a saturated set, cannot be 
  carried. Fortunately, the last item of Proposition \ref{diffeo} allows us to carry the property that the pair $(\F,\L)$
  is made of isotopic foliations. This is the reason why we introduce the notion  of \emph{weak isotopy}.
  
  \begin{defn}\label{weak}
  A \emph{weak isotopy} of $\L$ to $\F$ in $M$ is an alternating finite sequence of Dehn modifications---right or 
  left---and isotopies,
  $\M_1,\left(\psi_1^t\right), \ldots,\M_{k}, \left(\psi_k^t\right),\ t\in [0,1],$  and a sequence of pairs
  of foliations, $(\F,\L), (\F_1,\L_1), \ldots, (\F_k, \L_k)$, 
  fulfilling the following recursive conditions for every integer $j\in[1,k]$.
  \begin{enumerate}[{\rm (i)}]
  \item $\M_j$
  is a Dehn modification adapted to the pair of foliations $(\F_{j-1}, \psi^1_{j-1}(\L_{j-1}))$;\footnote{Here, $(\F_0,\L_0)=(\F,\L)$.}
  \item  $\M_j$ carries  $\F_{j-1}$ to $\F_j$ and $ \psi^1_{j-1}(\L_{j-1})$ to $\L_j$;
  \item $\F_k=\psi^1_k(\L_k )$.
  \end{enumerate}
  {\rm In this language, Theorem (A') reduces to:} $\F$ and $\L$ are weakly isotopic. \\
\end{defn}

  Now is the announced application of Dehn modification dealing with connecting orbits of saddles. 
  One should first mention that, as for a one-parameter family of functions on a Riemannian manifold, 
  generically on the metric the vector field $\nabla_\L f$ has finitely many
  orbits connecting two saddles; moreover there are connecting orbits neither from a max-saddle inflection to 
  a saddle nor from a saddle to a min-saddle inflection.
  We now include this property in the definition of \emph{excellence} and assume it
  in the remainder of  this paper. We recall that, by Proposition \ref{cancel-X}, we may assume that
   \emph{no saddles are of type $X$}.

  \begin{prop} \label{connectY-lamb} In this settting there exists a weak isotopy of $\ell$  to a new function $\ell'$ in $M$
  such that all of its $Y$-saddles are located above its $\la$-saddles  with respect to the order
  of their $f$-values.					
 \end{prop} 

\proof Without $X$-saddles the type of saddles is constant on every saddle arc.  The two ends 
of such an arc are inflections of the same Morse index:\footnote{For the index of an inflection 
see Subsection \ref{regul}.} 
if this index is $0$ (resp. $1$) the saddles are of type $\la$ (resp. $Y$)
as it can be seen near the inflection.

The matter is to destroy the $\L$-gradient connecting orbits $Y\to\la$. Indeed, assume the metric is chosen
so that there are no such 
connections. Consider  
a level $z_0$ above all $\la$-saddles and a contact arc $\al$ of $Y$-saddles.
By assumption on the metric,
the  ascending saturated set $S_{z_0}(\al)$ avoids every closure of $\la$-saddle arcs. 
Of course, $S_{z_0}(\al)$ could meet another $Y$-saddles not belonging to $\al$ 
or maxima and inflections of Morse index $1$.  
This does not matter, all thses contacts have to be 
incorporated to $S_{z_0}(\al)$.
An isotopy along $S_{z_0}(\al)$ realizes what is 
wished for $\al$, and so on for the other $Y$-saddle arcs.

To bypass the connecting orbits we use  Dehn modifications. At the beginning, we have finitely many 
connecting orbits $Y\to\la$ lying in distinct leaves $L_1, L_2,\ldots, L_n$ where each $L_i=\ell^{-1}(y_i)$
is a generic leaf. 
There exists $\ep>0$ such that $\hat L_i:= \ell^{-1}([y_i-\ep, y_i+\ep])$ is still made
of generic leaves.

The connecting orbit on  $L_i$ goes from $Y$-saddle $s_i$ to $\la$-saddle $s'_i$. Denote by $\al_i$ (resp. $\al'_i$)
the arc of contacts that contains $s_i$ (resp. $s'_i$). Set $\hat \al_i:= \hat L_i\cap \al_i$ and 
$\hat\al'_i:=\hat L_i\cap \al'_i$; all are named \emph{saddle intervals}.
 If $\ep$ is small enough there exists $\eta$ such that every interval in the collection
$\left(f(\hat\al_i), f(\hat\al'_j)\right)_{i,j\in[1, n]}$ has a length less than $\eta$ and every two of them
have a mutual distance larger than $\eta$. In addition, one may require, for every $j\in[1,n]$ and every pair 
$(a,b)$ of distinct contact points in a same leaf from $\hat L_j$,
\begin{equation}\label{eta}
\vert f(a)-f(b)\vert >2\eta.
\end{equation}
\begin{center}
 \begin{figure}[h]
\includegraphics[scale =.6]{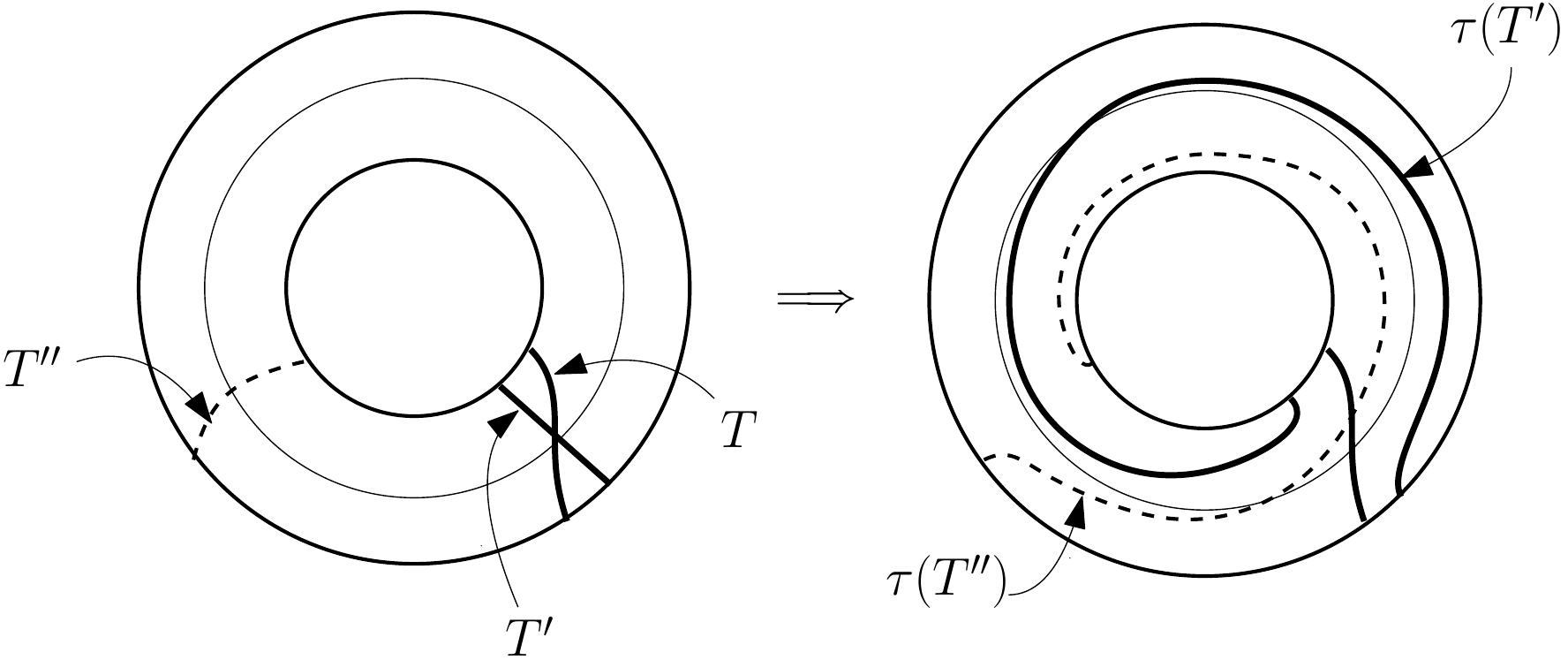}
\caption{\small $ T''$ is marked by $S_{z_1}(\hat\al'')$ where $\hat\al'' $ is a saddle interval of type $\la$ in $\hat L_1$ and $f(\hat\al'')>f(\hat\al'_1)$.
}\label{connecting-nb}
\end{figure}
\end{center}

Say $f(s_1)> f(s_i)$ for every $i>1$. Cut $\hat L_1$ at the level $z_1:=f(s'_1)-\eta$. By condition (\ref{eta}),
this level set of $\hat L_1$ is a finite collection of annuli. Among them, call $A_1$ the one which meets the 
connecting orbit $s_1\to s'_1$. 
The foliation $\L$ induces on $A_1$ a foliation by circles; the central circle lies in the leaf $L_1$.
Knowing that $s_1$ is of type $Y$ and $s'_1$ of type $\la$
one sees in $A_1$ two arcs $T$ and $T'$ joining the two boundary circles:
$T= A_1\cap S_{z_1}(\hat\al_1)$ and $T'= A_1\cap S_{z_1}(\hat\al'_1)$. Here, the first satured set is ascending
and the second one is descending. 
Moreover, $T\cap T'$ is the point where the connecting orbit crosses $A_1$---because 
the connecting orbits are on  isolated leaves.

 Applying  a Dehn twist  $\tau$ to $T'$---left or right depending on the sign of intersection $\langle T,
 T'\rangle$---together with the corresponding Dehn modification to $M$ kills the connecting orbit since
  $\tau(T')\cap T=\emptyset$.

However, new connecting orbits could have  been  created in the foliated domain $\hat L_1$ 
by this first Dehn modification (this is the case in Figure 
\ref{connecting-nb} due to  $\tau(T'')\cap T\neq\emptyset$.) 
One solves this new difficulty as follows.

Look at all saddle intervals of type $\la$ in $\hat L_1$ 
 at level higher than $z_1$. 
Name them $\hat\beta_1$, ..., $\hat\beta_k$ such that $f(\hat\beta_1)< f(\hat\beta_2)<\cdots<f(\hat\beta_k)$.\footnote{
 Inequalities of intervals mean that their are disjoint and ordered as it is written.}
The number of them is not 
affected by forthcoming Dehn modifications, it is kept invariant (Remark  \ref{important}.)

Say that $\hat \beta_1\subset \hat L_1$ is the lowest saddle interval of type $\la$ being $\L_\tau$-gradient 
connected to $\hat\al_1$; if there is no such a connecting orbit, consider $\hat\beta_2$, and so on. 
Take a proper annulus $A'_1 \subset\hat L_1$, at the  level $z'_1:=\inf f(\hat\beta_1)-\eta$, that crosses a 
connecting orbit from $\hat \alpha_1$ to $\hat\beta_1$. 
 
 The figure is similar to Figure \ref{connecting-nb}, but $S_{z'_1}(\hat\beta_1)$ marks a short arc $T'(\hat\beta_1)$ 
 joining the two
 boundary components of $A'_1$ while $S_{z'_1}(\hat\al_1)$ marks a long arc 
 like $\tau^{-1}(T)$. 
 These two arcs are transverse to the foliation $\L_\tau\cap A'_1$ and,
up to isotopy, there is only one intersection point between them. One can repeat the trick of Dehn modification to cancel
the connection $\hat\al_1\to \hat\beta_1$.
 Of course, new connecting orbits are possibly created,
 but only from $\hat\al_1$ to $\hat\beta_j$, $j>1$. So, the process goes on, without increasing the complexity.
 
 One can destroy all connecting orbits from $\hat\al_1$ to $\hat\beta_1\cup\cdots\cup \hat\beta_k$.
 As a result, after this long series of Dehn modifications, the ascending saturated set $S_{z_0}(\hat\al_1)$
 avoids every $\la$-saddle, where $z_0$ has been chosen at the very beginning of the proof. 
 If the saddles $(s_i)_{i=1}^n$, from which one connecting orbit goes up  in the leaf 
 $L_i$, are ranked in the decreasing order 
 of their $f$-values $\left(f(s_i)\right)_1^n$,  the process continues with the interval  $\hat\al_2$ and so on.\bull

 \section{Proof of Theorem (A')}\label{main-proof}

 \begin{rien}\label{ordered}{\rm
So far, we got that every saddle is of type $\la$ or $Y$---that is, no $X$-saddle---
 and, up to a weak isotopy applied to the function $\ell$,
 every $Y$-saddle is located at a level higher to all $\la$-saddles. Recall that an inflection adhering to a contact arc
 of $\la$-saddles (resp. $Y$-saddles) is of type saddle-min (resp. saddle-max). 
 Say that the $\la$-saddles 
 (resp. le $Y$--saddles) are
 located in $\{0<f<1/2\}$ (resp. $\{1/2<f<1\}$.)
 
 A way to achieve the proof of Theorem (A') is, by Lemma \ref{moser-trick2}, to kill every negative contact of the 
 pair $(f,\ell)$. This requires some algorithm since the leaves of $\L$ could have a very tricky topology
as embedded surface in $S^2\times[0,1]$.  In \cite{lauden}, this algorithm was based on words in an alphabet 
with four letters. Here, it is much simpler; it is based on two topological concepts: \emph{basin} and \emph{co-basin}
that we are going to define.
}
 \end{rien}

\begin{defn} Let $L=\{\ell=u_0\}$ be a leaf of $\L$ and let $m\in L$
be a negative contact of index $0$. The  \emph{basin} of $L$ determined by $m$ is the maximal 
3-ball $B(m)$ in $M=S^2\times [0,1]$ whose boundary is made of two parts: the first part is a horizontal disc 
$\p_\F B(m)$ whose boundary curve has exactly one negative $\la$-saddle if $L$ is a generic leaf or 
one of the following two situations: a negative inflection in $\p_\F B(m)$ or two saddles one of them being negative;
the second part 
is a disc $\p_\L B(m)$ in $L$ 
with $m$ as unique negative center and no other negative contact with $\F$ in the interior of $\p_\L B(m)$. 
Its \emph{threshold} consists of every 
negative saddle---generically unique---in the boundary curve of $\p_\L B(m)$. 
 \end{defn}
 
 In contrast, $\p_\L B(m)$ may have many
positive contacts. 
 If $z_0$ denotes the level of $\p_\F B(m)$ and $u_0$ denotes the $\ell$-value of the negative saddle(s), we have $B(m)\subset \{f\leq z_0\}\cap\{\ell\leq u_0\}$.
 
  \begin{defn}A basin $B(m)$ being given as above, a \emph{co-basin} joint to $B(m)$ is a maximal ball $B^*$
 in $\{f\leq z_0\}\cap\{\ell\geq u_0\}$ whose boundary is made of two parts:
 the first one is a horizontal disc $\p_\F B^*$; the second one is a disc $\p_\L B^*$ in the 
 $\L$-boundary of $ B(m)$. Its \emph{threshold} is the  unique\footnote{The uniqueness of its threshold holds 
 since there are no $X$-saddles (compare Figure \ref{saddle-nb}.)} $\la$-saddle which has  
 only one descending separatrix
 lying in $\p_\L B^*$ and the other in $\p_\L B(m)\smallsetminus B^*$.
 \end{defn}
 
 Somehow, a co-basin is like a pocket with respect to a basin. 
 Many co-basins could be glued to a basin.
 In our setting where the $Y$-saddles are above all $\la$-saddles, $\p_\L B^*$ cannot have maxima.
 
  \begin{center}
 \begin{figure}[h]
 \includegraphics[scale =.6]{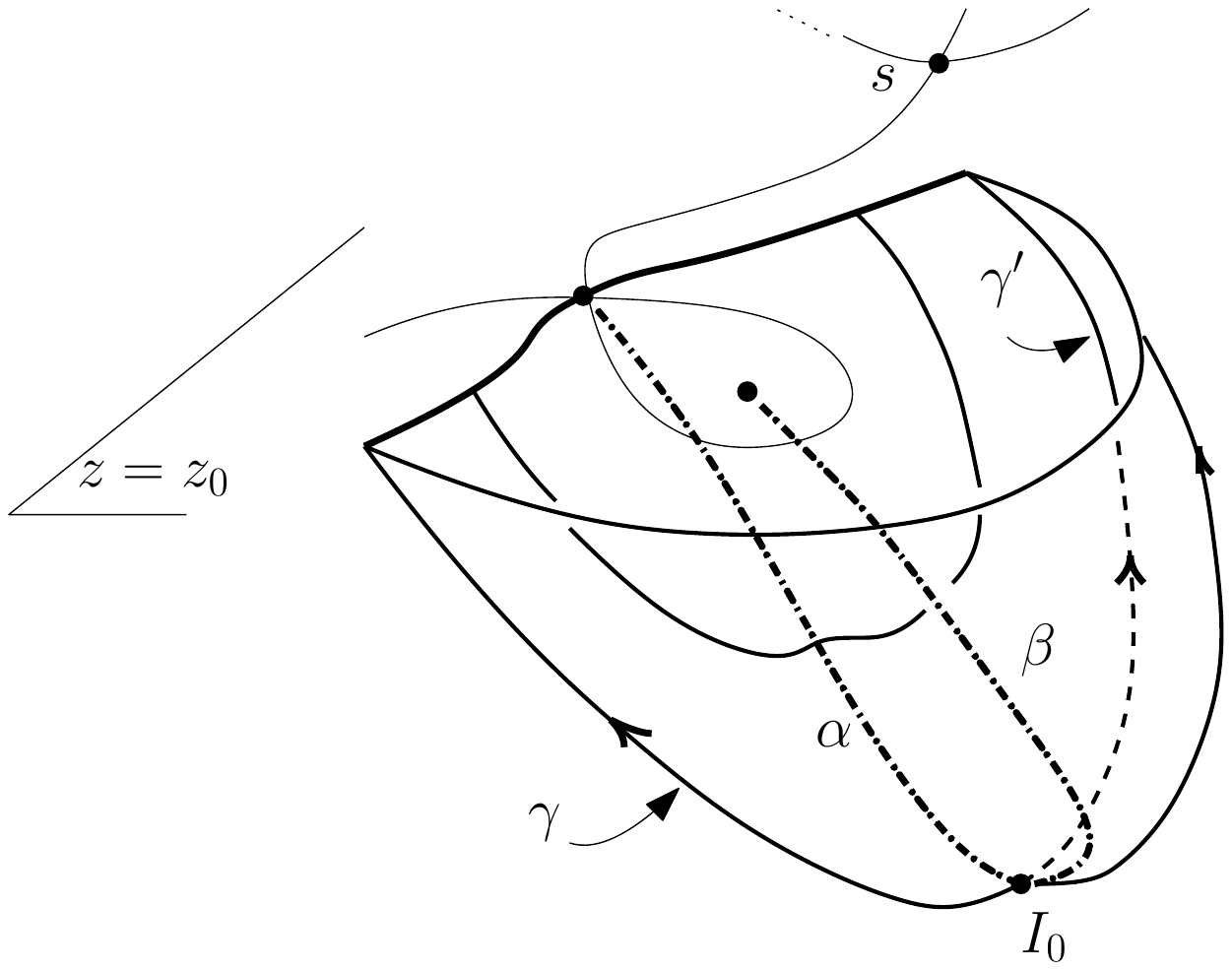}
\caption{ \small $I_0$ is a birth positive saddle-min inflection; $\al$ and $\beta$ are respectively the saddle 
and the min locus emanating from $I_0$; $\ga$ and $\ga'$ form the boundary of the unstable manifold of $I_0$.
The saddle $s$ is negative.}
\label{satur-nb}
 \end{figure}
 \end{center}


 \begin{rien}\label{complexity}{\sc Lower complexity.} {\rm Given a pair $(f,\ell)$ of functions without critical points,
 we define its \emph{lower complexity}
 $\kappa^-(f,\ell)$  as the pair $(\mu,\nu)$ of non-negative integers
 arranged in lexicographical order and defined in the following way.\footnote{This pair is meant to measure
 the topological complexity of $\L$ in the domain $\{0<f<1/2\}$.}
 
  The entry $\mu$ is the number of contact arcs of negative minima.
 
 Let $I_0$ be the \emph{upper}, unique by excellence, negative inflection which is the \emph{birth} of an  arc $\beta$ of negative minima.
Every minimum $m\in\beta$ defines a unique basin. 
  A minimum  $m\in \beta$, distinct from a cancellation inflection,
 is said to be \emph{accidental} if one of the following cases happens:
 \begin{itemize}
 \item[--] the boundary of $\p_\F B(m)$ contains either two saddles or a negative inflection---necessarily a cancellation 
 inflection---by definition of $I_0$,
  \item[--] $\p_\L B(m)$ contains a positive inflection. 
  \end{itemize}
  
 By definition, the entry $\nu$ of the lower complexity is the number of accidental minima in $\beta$.
  }
   \end{rien}
  
  If $\mu=0$ then there are neither negative minima nor negative $\la$-saddles anymore
  and, by convention, $\kappa^-(f,\ell)$ is equal to $(0,0)$; then the set of negative contacts is empty in $\{f<1/2\}$. 
   
 Similarly, some \emph{upper complexity}, $\kappa^+(f,\ell)$, is defined  that measures the topological complexity
 of $\L$ in the domain $\{1/2 <f<1\}$. Its vanishing means that there are no negative contacts in  $\{f>1/2\}$.
 By symmetry, any result about negative $\la$-saddles/minima applies to negative $Y$-saddles/maxima---and 
 conversely. Therefore we
are focusing on the first case. 

 \begin{lemme} \label{no-secondary}
 By the choice of $I_0$, denoting by $\beta$ the arc of minima that $I_0$ generates, 
 the following holds.
 \begin{enumerate}[\rm (1)]
 \item Let $B^*$ be a co-basin joint to a basin $B(m)$ with $m\in \beta$.
 Every (piece of) leaf in the interior of $B^*$ contains no negative minimum.
 \item The positive contact  arcs descending  from  contact points in the interior of $\p_\L B^*$ 
 descends to an inflection in the interior of $B(m)$. In particular, $m$ is the absolute minimum
 of the basin that it defines.
 \end{enumerate}
  \end{lemme}

\proof  
 1) Let $L$ be a leaf passing through $B^*$ and $m'$ be a negative minimum contact in $L\cap B^*$.
The contact arc $\beta'$ passing through $m'$ cannot cross $\p_\L B^*$ which only contains positive contacts. So, 
$\beta'$ descends to an inflection in $B^*$.
 But such an inflection should be higher than $m$, the negative minimum
of $B(m)$, that itself is higher than $I_0$. This contradicts the choice of $I_0$.

2) Suppose  such an arc $\ga$ crosses $\p_\L B(m)$ twice.  The integral of $d\ell$ on $\ga$, 
oriented in the direction of decreasing $z$, is positive. This prevents $\ga$ from having two distinct 
points on the same leaf.
Hence the lower point of $\ga$ is an inflection in the interior of $B(m)$. \bull

The topic of the next lemma is  to clean up the basins we are going to deal with in the proof of Theorem (A'), 
namely to 
make sure that no co-basins are joint to them. Again, 
$I_0$ is the \emph{upper} birth negative 
saddle-min inflection and $\beta$ is the contact arc of minima that $I_0$ generates. 

 \begin{lemme} \label{cleaning} In this setting, let $I'$ be  a \emph{positive} birth saddle-min inflection located in
 $\p_\L B(m)$ for some 
 $m\in \beta$. 
 Then
  there exists a \emph{contact conjugating}\footnote{A diffeomorphism $\Phi$ of $S^2\times[0,1]$
   is said to be  contact conjugating 
   if $\Phi$ maps every contact point of the pair $(\F,\L)$  to a contact point of the pair $(\F,\Phi(\L))$ and conversely.}
  isotopy
 $\Phi_t$, $t\in [0,1]$, fulfilling the following:
 \begin{enumerate}[{\rm (i)}]
 \item the isotopy is supported in $\{0<f<1/2\}$ and keeps the negative contacts fixed; 
 \item $\Phi_1(I')$ is away from every so-called $I_0$-basin, that is defined by some minimum in $\beta$. 
 \end{enumerate}
   After having iterated such isotopies, every $I_0$-basin is free of positive contact.
  \end{lemme}
  
  \proof For the first part we recall that when the negative threshold of a basin goes down on its own contact arc,
  meanwhile
  the positive threshold of the joint co-basin generated by $I'$ goes up on its own contact arc. Hence, there is a level 
  $z_{I'}$ which is the common level of both the threshold of some basin $B(m')$, $m'\in\beta$, 
  and the threshold 
  of the joint co-basin $B^*_{I'}$ generated by $I'$.
  
  So, it is natural to consider the ascending saturated set $\Si':=S_{z_{I'}+\ep}(I')$ for some small enough
   $\ep$. 
  Here, we recall that every ascending separatrix of a  $\la$-saddle reaches the level $\{z=1/2\}$ except 
  if it is bounded from above by a cancellation inflection or another $\la$-saddle---which 
  will be incorporated to $\Si'$. \smallskip
  
  
  \nd{\bf Claim.} {\it $\Si'$ does not approach the negative contacts.}\smallskip

  Indeed, by Lemma \ref{no-secondary}, there is no negative contact locus in the interior 
  of $\Si'$---what is  
  equivalent to the interior of a co-basin. What about $\al$, the locus of saddles emanating from $I_0$
  when the level is close to $z_{I'}$? Let $s'$ denote the threshold of the basin $B(m')$ and let $s'(z)$ 
  denote
  that nearby threshold at level $z$ close to $z_{I'}$. 
  Since the crossing can happen, certainly $s'(z)$
  has no connecting orbits with points of $\Si'$ if $\ep$ is small enough and $z$ ranges in 
  $(z(I')-\ep,z(I')+\ep)$.\footnote{Indeed, for $z=z_{I'}$, the threshold has a separatrix coming from $m'$
   for some choice of the Riemannian metric. Hence, the same holds in nearby basins.}
   So, for such a $z$, the negative saddle $s'(z)$ in not in the closure of $\Si'$. \bull
  
  By Lemma \ref{collapse-lem}, there is an isotopy along $\Si'$ that pushes it 
  above the level $z(I')$. Its support is located in a neighborhood of $\Si'$. 
  The time one of this isotopy fulfills the demand. \\
  
  About the last claim, it is sufficient to recall that there are only finitely many such positive birth inflections. 
  Let $I''$ be one of them.\footnote{Possibly, $I"$ belongs to $\Si'$. Hence, $\Phi_1(I'')$ is already away from any basin
  and hence this inflection can be skipped.}
   One can consider its saturated set $\Si''$ with respect to $\Phi_1(\L)$ up to a level which is 
  determined as previously by the crossing argument. An isotopy along $\Si''$ pushes $I"$ away from all basins
 and does not destroy what was gained in the first step. So, cumulating the isotopies related to each positive
 inflection that is initially contained in a basin makes all basins $B(m), \ m\in\beta$,
   free from positive contacts. The order in which the concerned inflections are numbered is irrelevant.
 \bull
 
 Now we are ready for the decisive part of the proof of Theorem (A'). It will consists of decreasing the
 lower complexity 
 $\kappa^-(f,\ell)$ until it vanishes. By the symmetry $z\to 1-z$, the same holds true about the upper  complexity.

 \begin{prop}\label{lambda2}  
 If the set of negative minima contacts of the pair $(f,\ell)$ is non-empty
 there exists an isotopy supported in $\{0<f<1/2\}$ that carries $\ell$ to a  function $\ell'$ such that 
 \begin{equation}\kappa^-(f,\ell')<\kappa^-(f,\ell).
 \end{equation}
 \end{prop}
 
 \proof There are several cases depending on how the complexity is topologically made. 
 
  1) Assume $\mu>0$ and $\nu=0$. We continue with
 the upper negative birth saddle-min inflection $I_0$ and the arc $\beta$ of index 0 that it generates.
 By Lemma \ref{cleaning}, every basin $B(m)$, $m\in\beta$, is free of positive contact.
 In this case, for every minimum $m\in \beta$ the cone $C_\L(m)$ is \break 
 standard.\footnote{See Subsection \ref{notation-ball} for notation.}
  
  If the negative $\la$-saddle arc $\al$ emanating from $I_0$ and the minimum arc $\beta$ close up in a common 
 cancellation inflection point then the union $\al\cup\beta$ forms a \emph{simple loop}  in the sense of Definition 
 \ref{def-simple}. By Proposition \ref{cancel-loop} there is an isotopy supported in a neighbourhood of 
 $\cup_{m\in \beta}C_\L(m)$ that eliminates this simple loop and decreases $\mu$ by $1$.
 
\begin{center}
 \begin{figure}[h]
 \includegraphics[scale =.6]{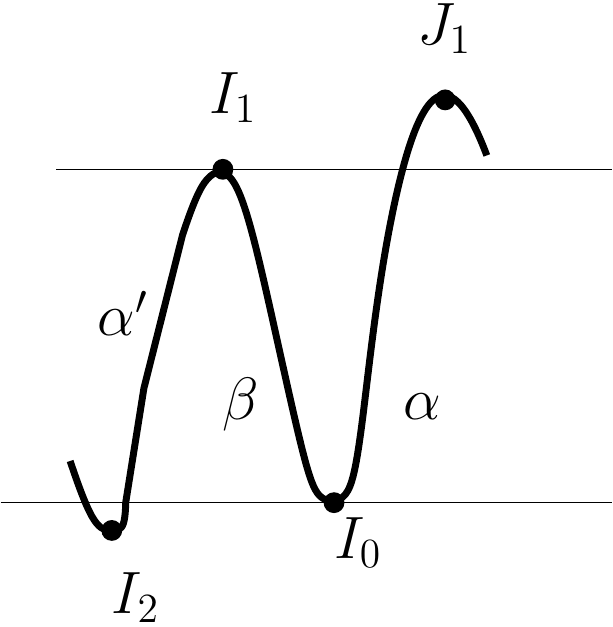}
\caption{\small The horizontal lines stand for  level sets.}\label{xxx0}
 \end{figure}
 \end{center}

 Let us show that there are no other configurations with $\nu=0$. Let $J_1$ be the inflection ending $\al$;
  let $I_1$, distinct from $ J_1$,
 be the inflection ending  $\beta$
 and let $I_2$ be the birth inflection of the $\la$-arc, named  $\al'$, descending from $I_1$. 
 We have $f(I_0)>f(I_2)$ by definition of $I_0$ (Figure \ref{xxx0}.) 
 Moreover,
 $f(J_1)>f(I_1)$; if not, the unique separatrix converging to
  $J_1$ starts from a minimum 
 $m_1\in \beta$, and hence $\nu>0$.
  
   When $\nu=0$ and $f(J_1)>f(I_1)$, 
every $m\in \beta$ ends one separatrix descending from some saddle $s(m)\in\al$
 and one separatrix descending from some saddle $s'(m)\in\al'$. The latter two claims
 prove similarly: the existence of $s(m)$ is clear near the inflection $I_0$ and extends to 
 $\beta$ since no accident happens along this arc; the existence of $s'(m)$ is clear near
 $I_1$ and also extends to $\beta$ for the same reason.
  

 When  $m$ is close to $I_0$ one has $f(s'(m))>f(s(m))$ and when $s'(m)$ is close to $I_1$ one has 
 $f(s'(m))<f(s(m))$. Then there is $m_0\in \beta$ such that 
 \begin{equation}
 f(s(m_0))=f(s'(m_0)). \label{eq-I_0}
 \end{equation} 
 Hence, $C_\L(m_0)$ is in the 
 \emph{saddle-center-saddle} configuration (Definition \ref{scs-conf}) 
  and $\nu$ is positive. Contradiction.

  2)  We assume $\nu>0$ in the rest of the proof. 
 Consider the first accidental minimum $m_0$ in $\beta$ when coming from 
 $I_0$ and assume that the curve $\p C_\L(m_0)$ contains two saddles $s_0$ and $s'_0$,
 necessarily negative by 
 the ``cleaning'' Lemma \ref{cleaning}.
 Then we have the saddle-center-saddle configuration; item (3) from Definition 
 \ref{scs-conf} is fulfilled by the choice of $I_0$ and the cleaning of the $I_0$-basins.
 
  \begin{center}
 \begin{figure}[h]
 \includegraphics[scale =.7]{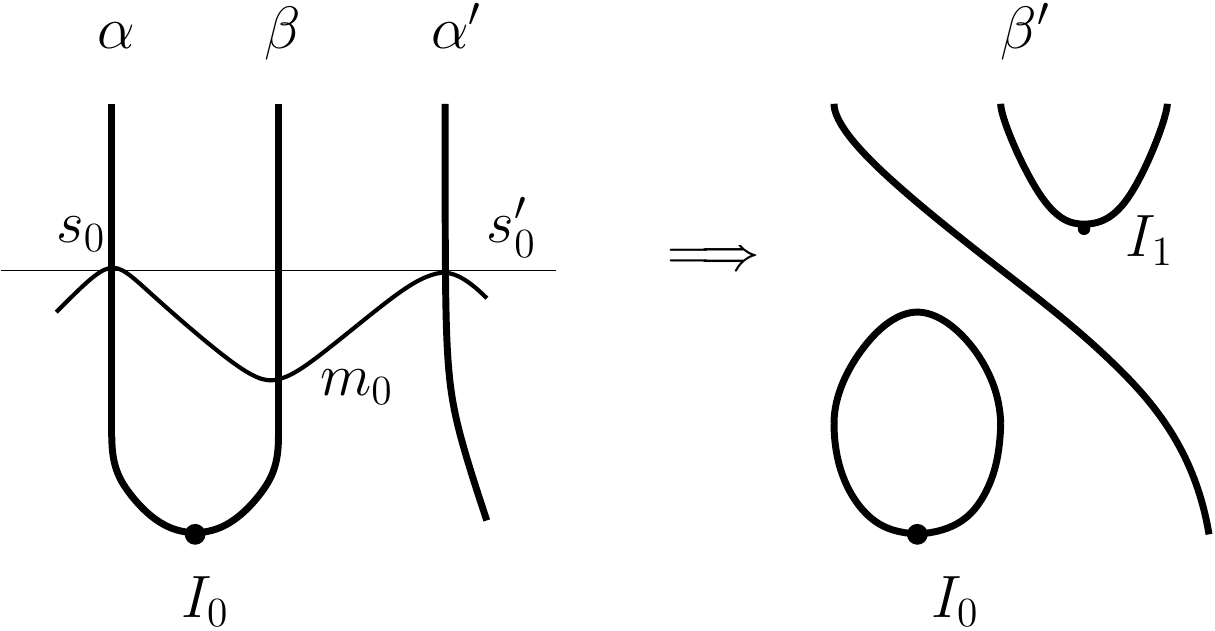}
\caption{Bypass creating a simple loop.}\label{bypass-cancel-nb}
 \end{figure}
 \end{center}
   Let $\al$ and $\al'$ denote the contact arcs containing $s_0$ and $s'_0$ respectively; let
 $s_t$, $s'_t$ and $m_t$ denote local parametrizations of $\al$, $\al'$ and $\beta$ respectively near $t=0$ such that 
 these three contact points are in the same leaf for every $t$ close to $0$.
 If $\al$ comes from $I_0$ and knowing that $m_0$ is the first accidental point on $\beta$, 
 one has $f(m_t)<f(s_t)<f(s'_t)$ for every $t<0$. 
 In other words, the $f$-values of $s_t$ and $s'_{t}$ are ordered as in Figure \ref{planar-fig-nb}.
 Thus every saddle in the sub-arc $(I_0, s_0)\subset \al$ is cancellable with the corresponding minimum in 
 $\beta$; for every small $t>0$, the pair $(s'_t, m_t)$  can be cancelled.
 Proposition \ref{scs-prop} applies and its effect is shown in Figure \ref{bypass-fig-nb}. In particular,
 this \emph{bypass} creates a simple loop which contains $I_0$ and a new birth inflection $I_1$ with the following 
 properties:
 \begin{itemize}
 \item $f(I_1)>f(I_0)$. 
 \item The minimum arc $\beta'$ emanating from $I_1$ has less accidental minima than $\beta$ since
 it coincides with $\beta$ in $\{f> f(I_1)+\ep\}$ for $\ep$ small enough (Figure \ref{bypass-cancel-nb}.)
 \end{itemize}

 Once the simple loop has been cancelled the number of birth inflections remains equal to the initial 
 $\mu$ but the entry $\nu$ of the 
 lower complexity decreases by $1$: the arc $\beta'$ has one less accidental minimum than $\beta$.\\

 \begin{center}
 \begin{figure}[h]
 \includegraphics[scale =.7]{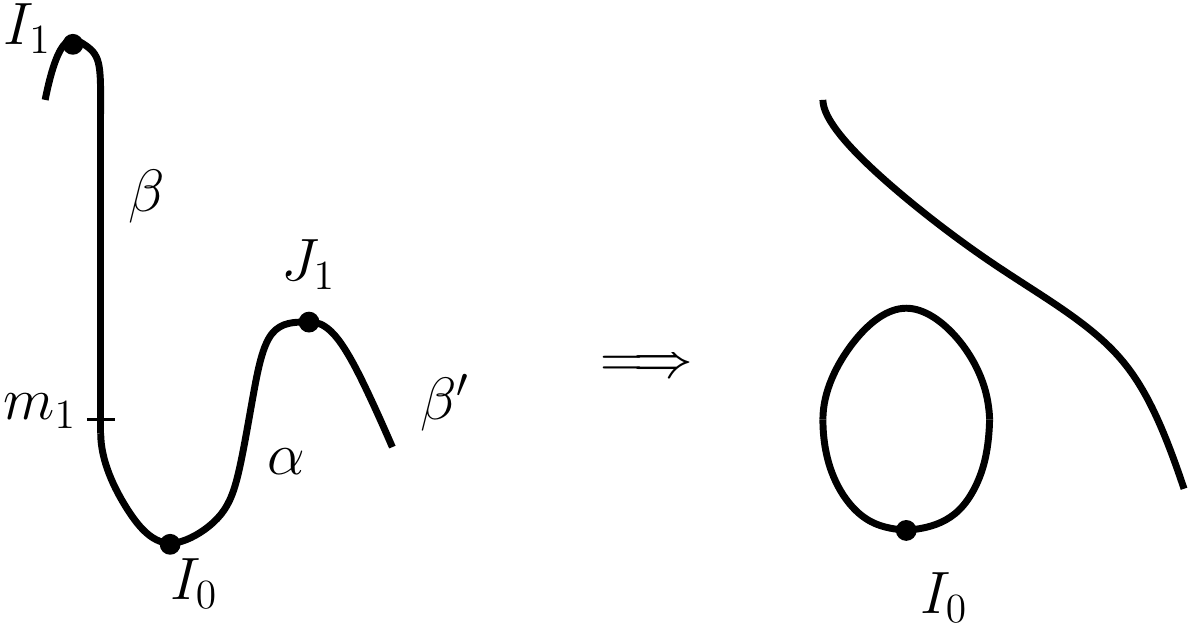}
\caption{Bypass in configuration min-saddle-min.}\label{bypass-cancel2-nb}
 \end{figure}
 \end{center}

 3) With the same notation, we assume $\nu>0$ and the first accidental minimum  $m_1\in \beta$ has a cone
 with an inflection point $J_1$ in its boundary. This inflection cannot be a birth since it lies at a higher level 
 than $I_0$. So, this
 is a cancellation inflection, which is a negative contact due to 
 the cleaning Lemma \ref{cleaning}.
 Let $\beta$ denote the index-0 arc starting from $I_0$ and let $I_1$ be its upper end. We have $I_1\neq J_1 $
 since an arc that is transverse to $\L$ cannot meet the same leaf twice.
 Here there are two cases.
 
 3-1) Assume that 
 the saddle arc $\al$
 starting from $I_0$ ends at $J_1$. Here, we are in the configuration
 min-saddle-min (Definition \ref{msm-conf}).  Let $\beta'$ be the arc of index 0 descending from $J_1$.
 The isotopy from Proposition \ref{msm-prop}
 modifies the contact arcs as, including their names, it is shown in Figure \ref{bypass2-nb}.
 The outcome is a simple loop containing $I_0$ and an arc of minima. 
 After the simple loop has been cancelled,
 the number of negative birth saddle-min inflections has decreased by one and hence the complexity decreases
 (Figure \ref{bypass-cancel2-nb}.)

 \begin{center}
 \begin{figure}[h]
\includegraphics[scale =.7]{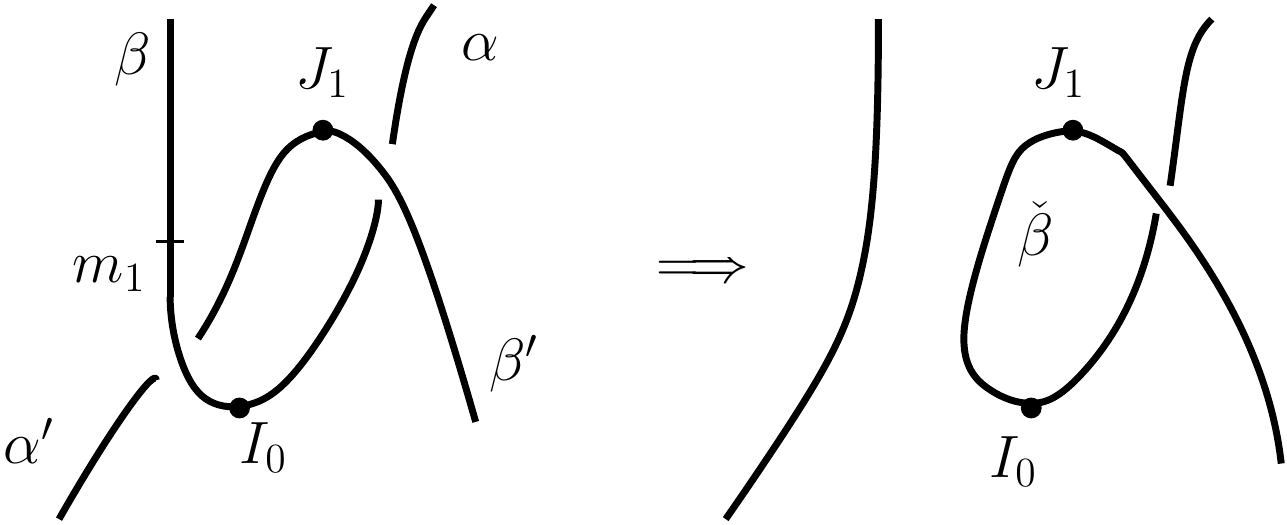}
\caption{Variant of Figure \ref{bypass-cancel2-nb}.}\label{bypass-cancel3-nb}
 \end{figure}
 \end{center}

 3-2) Assume that 
  $\al$ 
  does not end at $J_1$. Let $\al'$
 and $\beta'$ respectively denote the saddle arc and the minimum arc 
 ending at $J_1$. 
 The arcs $(\beta, \al',\beta')$ and the leaf $L$ that contains 
 $m_0\in \beta$ slightly below $m_1$ also
 presents a configuration min-saddle-min. The effect of the isotopy from Proposition \ref{msm-prop}
 does not produce a simple loop that could be cancelled. But it keeps $I_0$ as the upper negative birth inflection
 and the index-0 arc $\check\beta$ emanating from $I_0$ contains no accidental minimum.
 Therefore, the complexity decreases from $(\mu,\nu)$ to $(\mu, 0)$ (Figure \ref{bypass-cancel3-nb}.) \bull

 When the two complexities vanish, no negative contacts are remaining and
 the Moser trick (Lemma \ref{moser-trick2})  completes the proof of Theorem (A').\bull

\end{document}